# Some remarks and conjectures related to lattice paths in strips along the x-axis

Johann Cigler

Fakultät für Mathematik, Universität Wien

johann.cigler@univie.ac.at

**Abstract**

In the first part of this paper I give an elementary overview about some number sequences which count various sorts of lattice paths in strips along the x-axis and compute their generating functions in terms of Fibonacci and Lucas polynomials. In the second part I generalize these results by introducing suitable weights and study some special cases in more detail. In the course of this work I have been led to curious number triangles and various conjectures.

## 0. Introduction

Consider lattice paths in $\mathbb{Z}^2$ of length $n$ which start at the origin $(0,0)$ and have only up-steps $U:(i,j)\to(i+1,j+1)$ and down-steps $D:(i,j)\to(i+1,j-1)$. Equivalently consider (random) walks on $\mathbb{Z}$ which start at $0$, where at each step we go one unit up or down.

It is well known that the number of recurrent walks of length $2n$ is the same as the number of positive walks of length $2n$ or equivalently that the number of all lattice paths from $(0,0)$ to $(2n,0)$ is the same as the number of all non-negative lattice paths of length $2n$. The same holds for the number of all non-negative lattice paths of length $2n+1$ and the number of all lattice paths from $(0,0)$ to $(2n+1,-1)$.

Combining both results let us denote by $A_n$ the set of all lattice paths of length $n$ which start at the origin and end on height $0$ (if $n$ is even) or on height $-1$ (if $n$ is odd) and let $B_n$ be the set of all non-negative lattice paths of length $n$. Note that each path in $A_n$ has $\left\lfloor \frac{n}{2} \right\rfloor$ up- steps and $\left\lfloor \frac{n+1}{2} \right\rfloor$ down- steps. Then we have

$$|A_n| = |B_n| = \binom{n}{\left\lfloor \frac{n}{2} \right\rfloor}. \tag{0.1}$$

The main purpose of this paper is to study paths in $A_n$ with bounded heights.

Let $A_{n,k}$ be the set of all paths in $A_n$ which are contained in the strip $-\left\lfloor \frac{k+1}{2} \right\rfloor \leq y \leq \left\lfloor \frac{k}{2} \right\rfloor$ of width $\left\lfloor \frac{k}{2} \right\rfloor + \left\lfloor \frac{k+1}{2} \right\rfloor = k.$

For example $A_{5,3}$ consists of the following $8$ paths:

$DDUDU, DDUUD, DUDDU, DUDUD, DUUDD, UDDDU, UDDUD, UDUDD$

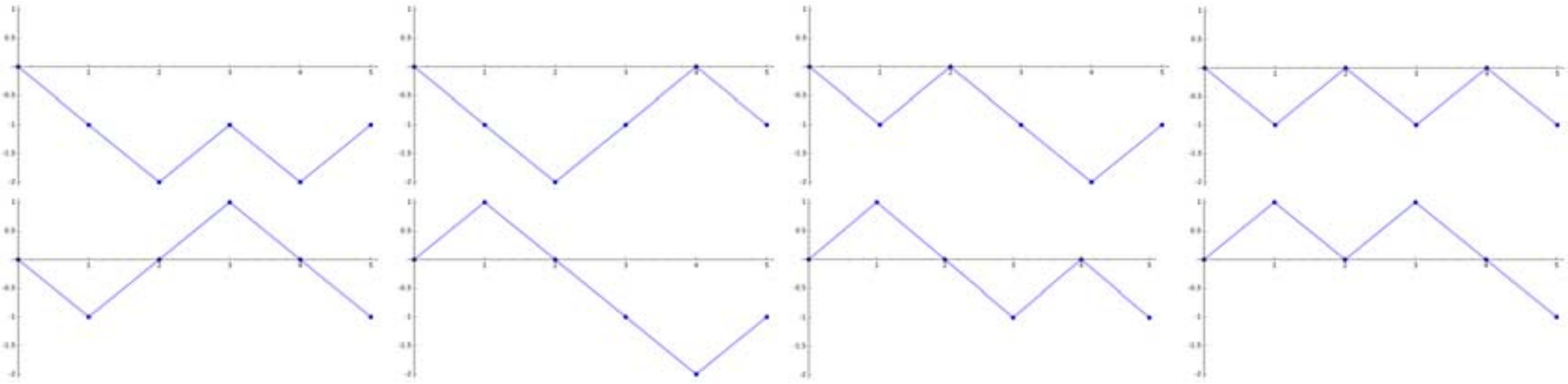

Figure 1

Let $B_{n,k}$ be the set of all paths in $B_n$ which remain in the strip $0 \leq y \leq k$.

For example $B_{5,3}$ consists of the following paths:

$UUUDU, UUUDD, UUDUU, UUDUD, UUDDU, UDUUU, UDUUD, UDUDU$

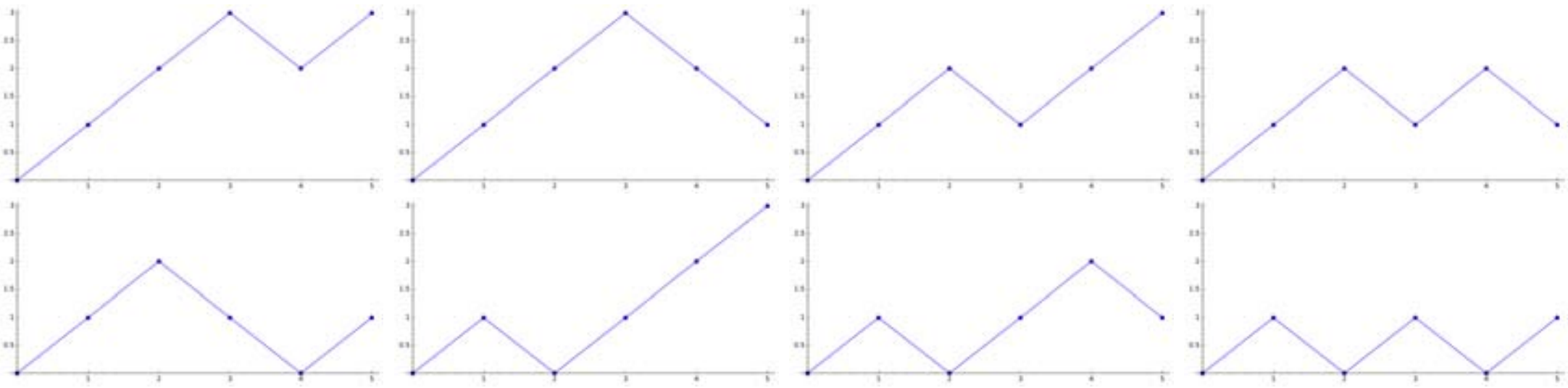

Figure 2

Then in analogy to (0.1) we have

$$\left|A_{n,k}\right| = \left|B_{n,k}\right| = \sum_{j\in\mathbb{Z}} (-1)^j \binom{n}{\left\lfloor \frac{n+(k+2)j}{2} \right\rfloor}. \tag{0.2}$$

This is essentially known, but we give another proof in the sequel.

The set $B_{n,k}$ can also be interpreted as the set of all walks on the path graph with vertices $\{0,1,\cdots,k\}$.

Let $M_k = \left(m(0,i,j,k)\right)_{i,j=1}^{k} = \left([|i-j|=1]\right)_{i,j=0}^{k}$ be the corresponding adjacency matrix. Note that $m(0,i,j,k)=1$ if $\{i,j\}$ is an edge and $m(0,i,j,k)=0$ else, i.e. $m(0,0,1,k)=1$, $m(0,k,k-1,k)=1$, and $m\left(i,i\pm 1,k\right)=1$ for $0<i<k$. All other entries are $0$.

For example

$$M_3 = \begin{pmatrix} 0 & 1 & 0 & 0 \\ 1 & 0 & 1 & 0 \\ 0 & 1 & 0 & 1 \\ 0 & 0 & 1 & 0 \end{pmatrix}.$$

If we denote by $m(n,i,j,k)$ the number of non-negative paths from $(0,i)$ to $(n,j)$ in $B_{n,k}$ then by the definition of matrix multiplication we get

$$M_k^n = \left(m(n,i,j,k)\right)_{i,j=0}^{k}. \tag{0.3}$$

For example

$$M_3^5 = \begin{pmatrix} 0 & 5 & 0 & 3 \\ 5 & 0 & 8 & 0 \\ 0 & 8 & 0 & 5 \\ 3 & 0 & 5 & 0 \end{pmatrix}.$$

In the first row we see that there are $5$ paths from $(0,0)$ to $(5,1)$ and $3$ paths from $(0,0)$ to $(5,3)$.

In the general case we get $\left|B_{n,k}\right| = \sum_{j=0}^{k} m(n,0,j,k)$.

A little thought also gives that

$$\begin{aligned} \left|A_{2n,k}\right| &= m\left(2n, \left\lfloor \frac{k+1}{2} \right\rfloor, \left\lfloor \frac{k+1}{2} \right\rfloor, k\right) \\ \left|A_{2n+1,k}\right| &= m\left(2n+1, \left\lfloor \frac{k+1}{2} \right\rfloor, \left\lfloor \frac{k-1}{2} \right\rfloor, k\right). \end{aligned} \tag{0.4}$$

To see this renumber the rows and columns from $-\left\lfloor \frac{k+1}{2} \right\rfloor$ to $\left\lfloor \frac{k}{2} \right\rfloor$.

Then the paths from $B_{2n,k}$ which start from $\left(\left\lfloor \frac{k+1}{2} \right\rfloor, 0\right)$ and end in $\left(2n, \left\lfloor \frac{k+1}{2} \right\rfloor\right)$ are mapped to the paths from $(0,0)$ to $(2n,0)$ in $A_{2n,k}$ and the paths from $B_{2n+1,k}$ which start from $\left(\left\lfloor \frac{k+1}{2} \right\rfloor, 0\right)$ and end in $\left(2n+1, \left\lfloor \frac{k-1}{2} \right\rfloor\right)$ are mapped to the paths from $(0,0)$ to $(2n+1,-1)$ in $A_{2n+1,k}$.

In our example $\left|A_{5,3}\right| = m(5,2,1,3) = 8.$

In the first part of this paper I give an overview about these numbers. Most of these results are known but perhaps my point of view gives a novel approach.

In the second part we consider the following weights instead of the numbers $\left|A_{n,k}\right|$.

Define a peak as a vertex preceded by an up-step $U$ and followed by a down-step $D$, and a valley as a vertex preceded by a down- step $D$ and followed by an up-step $U$. The height of a vertex is its $y-$ coordinate. The peaks with a height at least $1$ and the valleys with height at most $-2$ are called extremal points.

Thus for example the path *UDDUDDU*

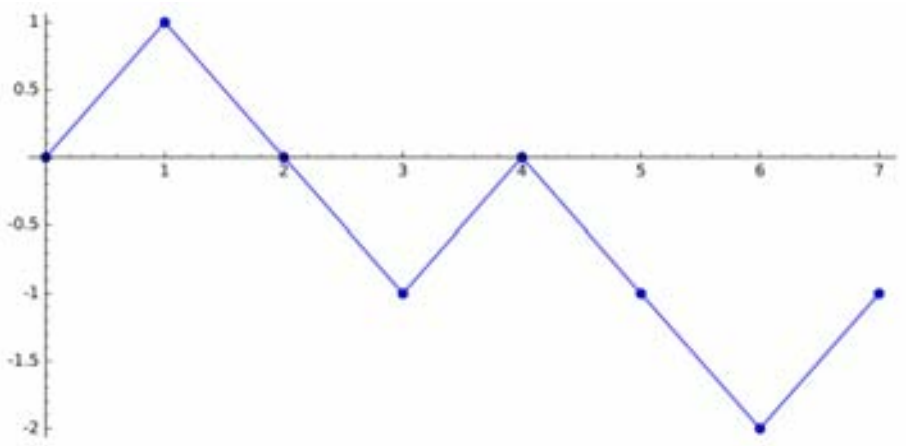

has two extremal points $\{(1,1),(6,-2)\}$.

Let $E(v)$ be the set of $x-$coordinates of the extremal points of the path $v$, $e(v) = \left|E(v)\right|$ the number of extremal points of $v$ and $\iota(v) = \sum_{i \in E(v)} i$ the sum of the $x-$coordinates of the extremal points. In [7] and [8] we defined the weight of $v$ by $w(v,q) = q^{\iota(v)} t^{e(v)}$ and considered the polynomials $w\left(A_{n,k}, q\right) = \sum_{v \in A_{n,k}} q^{\iota(v)} t^{e(v)}$. These polynomials are intimately connected with Rogers-Ramanujan type theorems.

In the present paper we consider only the case $q = 1$ and study the polynomials

$$a(n,k,t) = \sum_{v \in A_{n,k}} t^{e(v)} \tag{0.5}$$

in more detail.

In our example $A_{5,3}$ we get $a(5,3,t) = t^2 + t + t + 1 + t + t^2 + t + t^2 = 1 + 4t + 3t^2$.

The sequences $\left(a(n,k,t)\right)_{n\geq 0}$ satisfy recurrences of an order $\leq k$ and their generating functions are rational functions of Fibonacci and Lucas polynomials. For $n\to\infty$ they converge to $a(n,t)=\sum_{\ell\geq 0}\binom{\left\lfloor \frac{n}{2}\right\rfloor}{\ell}\binom{\left\lfloor \frac{n+1}{2}\right\rfloor}{\ell}t^{\ell}$. Therefore their Hankel determinants

$d(n,k,t)=\det\left(a(i+j,k,t)\right)_{i,j=0}^{n}$ vanish for $n\geq k$ and converge for each fixed $n$ to

$d(n,t)=\det\left(a(i+j,t)\right)_{i,j=0}^{n}=t^{\left\lfloor \frac{(n+1)^2}{4}\right\rfloor}$. Curiously the polynomials $d(n,k,t)$ for even $k$ are multiples of Hankel determinants of $a(n,t)$ and for odd $k$ multiples of Hankel determinants of Narayana polynomials $N_n(t)=\sum_{j=0}^{N}\binom{n}{j}\binom{n}{j-1}\frac{1}{n}t^{j}$.

As far as I know these polynomials have not been considered in the literature. The detailed study of some special cases led to many curious conjectures. Some have been found with the help of the Mathematica package Guess [15] by Manuel Kauers. Of great use has been The On-Line Encyclopedia of Integer Sequences OEIS [19]. If some results are already known please let me know that I can give due credit.

## 1. Background material

### 1.1. Lattice paths in strips along the x-axis

Let $A_{n,k}$ be the set of all lattice paths of length $n$ which start at $(0,0)$, stop on heights $0$ or $-1$ and are contained in the strip $-\left\lfloor \frac{k+1}{2}\right\rfloor \leq y\leq \left\lfloor \frac{k}{2}\right\rfloor$ of width $\left\lfloor \frac{k}{2}\right\rfloor+\left\lfloor \frac{k+1}{2}\right\rfloor=k.$

For $n\leq k$ all $\binom{n}{\left\lfloor \frac{n}{2}\right\rfloor}$ paths of length $n$ belong to $A_{n,k}$. Note that for odd $k$ the strips are not symmetric about the $x-$ axis.

By inclusion – exclusion it has been shown (see e.g. [7],[8] or [9] ) that

$$a(n,k):=\left|A_{n,k}\right|=\sum_{j\in\mathbb{Z}}(-1)^{j}\binom{n}{\left\lfloor \frac{n+(k+2)j}{2}\right\rfloor}. \qquad (1.1)$$

We will give a new proof of this result by showing that $\left|B_{n,k}\right|$ satisfies (1.1) in Proposition 1.4 and that $\left|A_{n,k}\right|=\left|B_{n,k}\right|$ in Corollary 1.3.

The set $A_{n,0}$ is empty for $n > 0$ which gives

$$\left|A_{n,0}\right| = \sum_{j\in\mathbb{Z}} (-1)^j \binom{n}{\left\lfloor \frac{n}{2} \right\rfloor + j} = [n = 0]. \tag{1.2}$$

For $k = 1$ the sets $A_{n,1}$ consist only of one path. If we denote a path by the sequence of its successive heights this unique path is $(0, -1, 0, -1, \cdots)$. We can write it also as $DUDU \cdots DUD$.

Therefore we have

$$\left|A_{n,1}\right| = \sum_{j\in\mathbb{Z}} (-1)^j \binom{n}{\left\lfloor \frac{n+3j}{2} \right\rfloor} = 1. \tag{1.3}$$

For example we have $\left|A_{7,1}\right| = \binom{7}{0} - \binom{7}{2} + \binom{7}{3} - \binom{7}{5} + \binom{7}{6} = 1 - 21 + 35 - 21 + 7 = 1.$

It is also clear that $B_{n,1}$ consists only of the path $UDUDU \cdots$.

The sets $A_{n,2}$ are $\{(0)\}, \{(0,-1)\}, \{(0,1,0),(0,-1,0)\}, \{(0,1,0,-1),(0,-1,0,-1)\}, \cdots$ or $\{\varepsilon, \{D\}, \{UD, DU\}, \{UDD, DUD\}, \{UDUD, UDDU, DUUD, DUDU\}, \cdots\}$ if we denote by $\varepsilon$ the trivial path, which gives by induction

$$\left|A_{n,2}\right| = \sum_{j\in\mathbb{Z}} (-1)^j \binom{n}{\left\lfloor \frac{n+4j}{2} \right\rfloor} = 2^{\left\lfloor \frac{n}{2} \right\rfloor}. \tag{1.4}$$

There is an easy bijection $\varphi$ between $A_{n,2}$ and $B_{n,2}$. Define $\varphi(D) = U$, $\varphi(UD) = UD, \varphi(DU) = UU$ and $\varphi(UDv) = UD\varphi(v)$ and $\varphi(DUv) = UU\,(\varphi(v))^*$, where $(\varphi(v))^*$ exchanges $U$ and $D$ in $\varphi(v)$.

Thus for example $\varphi(UDD) = UD\varphi(D) = UDU$, $\varphi(DUD) = UU\varphi(D)^* = UUD$, $\varphi(UDUD) = UDUD$, $\varphi(UDDU) = UD\varphi(DU) = UDUU$, $\varphi(DUUD) = UU\varphi(UD)^* = UUDU$, $\varphi(DUDU) = UU\varphi(DU)^* = UUDD$.

A very interesting case occurs for $k=3$. In this case we have $\left|A_{n,3}\right| = F_{n+1}$, where $F_n$ is a Fibonacci number.

The Fibonacci numbers $F_n = \sum_{k=0}^{\left\lfloor \frac{n-1}{2} \right\rfloor} \binom{n-1-k}{k}$ satisfy $F_n = F_{n-1} + F_{n-2}$ with initial values $F_0 = 0$ and $F_1 = 1$. The first terms are $0,1,1,2,3,5,8,13,21,34,\cdots$ (cf. OEIS [19], A000045).

Thus

$$\left|A_{n,3}\right| = \sum_{j\in\mathbb{Z}} (-1)^j \binom{n}{\left\lfloor \frac{n+5j}{2} \right\rfloor} = F_{n+1} = \sum_{k=0}^{\left\lfloor \frac{n}{2} \right\rfloor} \binom{n-k}{k}. \tag{1.5}$$

Since $A_{0,3} = \{0\}$ and $A_{1,3} = \{(0,-1)\}$ we see that the initial values are $F_1 = F_2 = 1$.

Consider now a path in $A_{n,3}$. Let the path be given by the sequence of its $y$ – coordinates. If the path ends with $(-1,0)$ or $(0,-1)$ then the path is the unique continuation of a path in $A_{n-1,3}$. If the path ends with $(0,1,0)$ or $(-1,-2,-1)$ then it is the unique continuation of a path in $A_{n-2,3}$. Therefore we have $\left|A_{n,3}\right| = \left|A_{n-1,3}\right| + \left|A_{n-2,3}\right|$. This together with the initial values gives $\left|A_{n,3}\right| = F_{n+1}$.

**Remark 1.1**

The formula $\sum_{j\in\mathbb{Z}} (-1)^j \binom{n}{\left\lfloor \frac{n+5j}{2} \right\rfloor} = F_{n+1}$ has been obtained by G.E. Andrews [1], but already in 1917 I. Schur [25] has studied the right-hand side of the identity

$$\sum_{k=0}^{\left\lfloor \frac{n}{2} \right\rfloor} q^{k^2} \begin{bmatrix} n-k \\ k \end{bmatrix}_q = \sum_{j\in\mathbb{Z}} (-1)^j q^{\frac{j(5j-1)}{2}} \begin{bmatrix} n \\ \left\lfloor \frac{n+5j}{2} \right\rfloor \end{bmatrix}_q \tag{1.6}$$

for $|q| < 1$. Here $\begin{bmatrix} n \\ k \end{bmatrix} = \begin{bmatrix} n \\ k \end{bmatrix}_q$ denotes a $q-$ binomial coefficient defined by

$\begin{bmatrix} n \\ k \end{bmatrix}_q = \dfrac{(1-q)(1-q^2)\cdots(1-q^{n-k})}{(1-q)(1-q^2)\cdots(1-q^k)}$ for $0 \le k \le n$ and $0$ else.

It is clear that (1.6) converges to (1.5) for $q \to 1.$ Therefore (1.6) is called a $q-$ analogue of (1.5).

If we let $n \to \infty$ in (1.6) we get the famous (first) Rogers-Ramanujan identity

$$\sum_{k\geq 0} \frac{q^{k^2}}{(1-q)(1-q^2)\cdots(1-q^k)} = \frac{1}{\prod_{j=1}^{\infty}(1-q^j)} \sum_{j\in\mathbb{Z}} (-1)^j q^{\frac{j(5-j)}{2}}. \tag{1.7}$$

For $k=3$ Thomas Prellberg [22] has found a simple bijection from $A_{n,3}$ to $B_{n,3}$:

First exchange $U$ and $D$, remove the last step and append the remaining path to the new first step $U$.

Thus for example we get $DDUDU \to UUDUD \to UUDU \to UUUDU$ or $UDDUD \to DUUDU \to DUUD \to UDUUD.$

From Figure 1 we get in this way Figure 2.

This map is obviously invertible. For a path in $B_{n,3}$ remove the first step $U.$ Then exchange $U$ and $D.$ This gives a path with heights between $-2$ and $1.$ Now append the uniquely determined step such that the path ends on height $0$ or $-1.$

Another bijection has been given by Helmut Prodinger [23].

It would be interesting to find simple bijections between $A_{n,k}$ and $B_{n,k}$ for $k>3.$

Some of the sequences $\left(a(n,k)\right)_{n\geq 0}$ occur in the literature in other contexts.

For small values of $k$ useful information can be found in OEIS [19]:

$\left(a(n,2)\right)_{n\geq 0}$ is A016116,

$\left(a(n,3)\right)_{n\geq 0}$ is A000045,

$\left(a(n,4)\right)_{n\geq 0} = \left(1,1,2,3,6,9,18,27,\cdots\right)$ is A182522,

$\left(a(n,5)\right)_{n\geq 0} = \left(1,1,2,3,6,10,19,33,\cdots\right)$ is A028495,

$\left(a(n,6)\right)_{n\geq 0} = \left(1,1,2,3,6,10,20,34,68,\cdots\right)$ is A030436,
$\left(a(n,7)\right)_{n\geq 0} = \left(1,1,2,3,6,10,20,35,69,\cdots\right)$ is A061551 and
$\left(a(n,8)\right)_{n\geq 0} = \left(1,1,2,3,6,10,20,35,70,125,\cdots\right)$ is A178381.

My interest in these topics has been aroused by the curious formula (1.5) for the Fibonacci numbers. In [6] I tried to put this identity into a general context in order to find an "explanation" of this formula. Among other things I proved that sums of the form (1.1) satisfy some simple recurrences. In the papers [10] – [12] I found simpler proofs and determined the generating functions of these numbers. From [20] I learned the interpretation of these sums as numbers of lattice paths which led to papers [7] – [9]. From this point of view formula (1.5) appears as a special case of the principle of inclusion – exclusion. Later I found in OEIS [19] for the above mentioned special cases of sequences $(a(n,k))$ the interpretation as walks in path graphs from which I got new insight into the situation. It turned out that some results which I had previously obtained were already known in other contexts. In the following pages I give an account of my present knowledge of this topic. Remarks and hints to the literature are very welcome.

## 1.2. Some other combinatorial models

**Proposition 1.1**

*The number $a(n,2k+1)$ counts all non-negative lattice paths starting from $(0,0)$ and ending in $(n,0)$, where besides up-steps and down-steps also horizontal moves $(i,0) \to (i+1,0)$ on height $0$ are allowed and the maximal height of a path is $k$.*

It is easy to find a bijection between these two lattice path models. Starting from the first model we map each up-step $(i,-1) \to (i+1,0)$ and each down-step $(i,0) \to (i+1,-1)$ into a horizontal move $(i,0) \to (i+1,0)$. The non-negative parts of the path remain unaltered and the negative paths $(i,-1) \to (j,-1)$ are reflected on the line $y = -\frac{1}{2}$ into a non-negative path $(i,0) \to (j,0)$. This map obviously has a unique inverse.

For example the set $A_{5,3}$ is transformed to

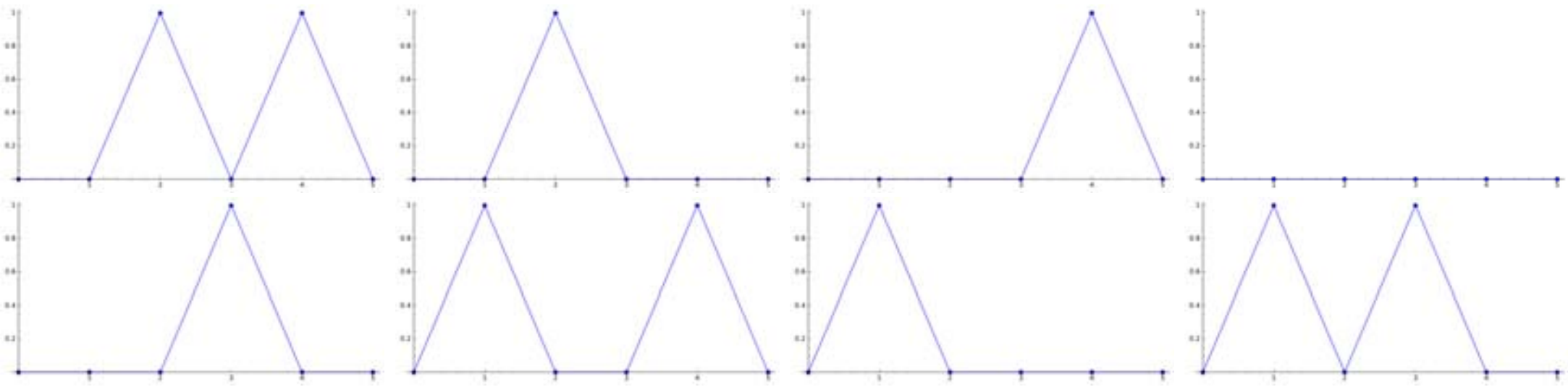

Figure 3

In this case the transformation has already been considered in [20].

**Corollary 1.1**

*Consider the graph with vertices $\{0,1,\cdots,k\}$ which arises by adjoining a loop at $0$ to the path graph with these vertices and let $R_k$ be its adjacency matrix. Then $a(n,2k+1)$ counts all walks of length $n$ which start and end in $0$. Therefore $a(n,2k+1)$ is the upper-left-most entry in $R_k^n$.*

Consider for example $2k+1=7.$ Then we get

$$R_3=\begin{pmatrix}1&1&0&0\\1&0&1&0\\0&1&0&1\\0&0&1&0\end{pmatrix} \text{ and e.g. } R_3^8=\begin{pmatrix}69&55&48&20\\55&62&27&28\\48&27&42&7\\20&28&7&14\end{pmatrix}.$$

Thus $a(8,7)=69.$

**Remark 1.2**

If we let $k\to\infty$ in Proposition 1.1, i.e. consider non-negative lattice paths starting from $(0,0)$ where besides up-steps and down-steps also horizontal steps $(i,0)\to(i+1,0)$ on height $0$ are allowed then the numbers $c(n,j)$ of all such paths ending in $(n,j)$ satisfy

$c(n,0)=c(n-1,0)+c(n-1,1)$ and $c(n,j)=c(n,j-1)+c(n,j+1)$ for $j>0.$

We get the following table (cf. OEIS A061554)

$$\left(c(n,j)\right)_{n,j\geq 0}=\begin{pmatrix}1&&&&&\\1&1&&&&\\2&1&1&&&\\3&3&1&1&&\\6&4&4&1&1&\\10&10&5&5&1&1\end{pmatrix} \tag{1.8}$$

.

It is easily verified that $c(n,j)=\binom{n}{\left\lfloor\frac{n-j}{2}\right\rfloor}$. Of course we have $c(n,0)=\binom{n}{\left\lfloor\frac{n}{2}\right\rfloor}$.

If we make the further assumption that $c(n,k+1)=0$ then we get by Proposition 1.1 that $c(n,0)=a(n,2k+1).$

**Proposition 1.2**

*The number $a(n,2k)$ counts all non-negative lattice paths starting from $(0,0)$ and ending in $(n,n\bmod 2)$, where the down-steps $(i,1)\to(i+1,0)$ are counted twice and the maximal height of a path is $k$.*

**Proof**

Each path has the form $P_1P_2\cdots P_jR$ where each $P_i$ is either non-negative or non-positive and starts and ends in height $0$ and $R$ is empty if $n\equiv 0 \bmod 2$ and negative if $n\equiv 1 \bmod 2.$ We reflect each non-positive $P_i$ on the $x-$axis and let the non-negative ones fixed. Thus each $P_i$ occurs twice. This is accounted for by counting the last down-step twice. The reflected path $R'$ ends in $(n,1).$

Thus for $k=4$ the numbers $\left(a(n,4)\right)_{n\geq 1}=\left(1,2,3,6,9,\cdots\right)$ count the paths

$\{U\}$, $\{U\boldsymbol{D}\}$, $\{U\boldsymbol{D}U,UUD\}$, $\{U\boldsymbol{D}U\boldsymbol{D},UUD\boldsymbol{D}\}$, $\{U\boldsymbol{D}U\boldsymbol{D}U,U\boldsymbol{D}UUD,UD\boldsymbol{D}U,UUDUD\},\cdots,$

where the bold letters indicate the steps which are counted twice.

If we let $k\to\infty$ and denote by $b(n,j)$ the weighted number of paths from $(0,0)$ to $(n,j),$ then we get $b(n,j)=\binom{n}{\frac{n-j}{2}}$ if $n-j$ is even and $b(n,j)=0$ else. This is OEIS A108044.

**Corollary 1.2**

*Consider the graph with vertices $\{0,1,\cdots,k\}$ which arises from the path graph by adjoining a second edge $(1,0)$ and let $S_k$ be its adjacency matrix. Then $a(2n,2k)=\left(S_k^{2n}\right)_{0,0}$ is the upper-left-most entry in $S_k^{2n}$ and $a(2n+1,2k)=\left(S_k^{2n+1}\right)_{0,1}=\left(S_k^{2n}\right)_{1,1}$. The last identity is clear because each path begins with an up-step.*

Consider for example $k=3.$

Then $S_3=\begin{pmatrix}0&1&0&0\\2&0&1&0\\0&1&0&1\\0&0&1&0\end{pmatrix}$ and $S_3^6=\begin{pmatrix}20&0&14&0\\0&34&0&14\\28&0&20&0\\0&14&0&6\end{pmatrix}.$

Thus $a(6,6)=20$ and $a(7,6)=34.$

**Proposition 1.3**

*Let $M_k=\left(m(0,i,j,k)\right)_{i,j=1}^{k}=\left([|i-j|=1]\right)_{i,j=0}^{k}$. Then $M_k^n=\left(m(n,i,j,k)\right)_{i,j=0}^{k}$, where $m(n,i,j,k)$ denotes the number of paths from $(0,i)$ to $(n,j).$*

*If we extend $m(n,i,j,k)$ by setting*

$$\begin{aligned}m(n,0,-2-j,k)&=-m(n,0,j,k),\\ m(n,0,j+k+2,k)&=-m(n,0,k+1-j,k),\end{aligned}\tag{1.9}$$

*then we get*

$$m(n,i,j,k)=\sum_{\ell=0}^{i} m(n,0,j-i+2\ell,k) \tag{1.10}$$

**Proof**

We know that $m(n+1,i,j,k)=m(n,i,j-1,k)+m(n,i,j+1,k)$ for $0\le i,j\le k$. Since there are no paths to $-1$ and to $k+1$ we also have $m(n,i,-1,k)=m(n,i,k+1,k)=0$. There is a uniquely determined extension of $m(n,i,j,k)$ for $j\in\mathbb{Z}$ which satisfies $m(n,i,j,k)=m(n-1,i,j-1,k)+m(n-1,i,j+1,k)$ for all $j\in\mathbb{Z}$. This satisfies (1.9).

To prove the proposition we show by induction that

$r(n,i,j,k)=\sum_{\ell=0}^{i} m(n,0,j-i+2\ell,k)$ equals $m(n,i,j,k)$.

For $n=0$ we have $r(0,i,j,k)=m(0,i,j,k)$. If we already know that $r(n-1,i,j,k)=m(n-1,i,j,k)$ then

$$r(n,i,j,k)=\sum_{\ell=0}^{i} m(n,0,j-i+2\ell,k)=\sum_{\ell=0}^{i}\left(m(n-1,0,j-i+2\ell-1,k)+\left(m(n-1,0,j-i+2\ell+1,k)\right)\right)$$
$$=r(n-1,i,j-1,k)+r(n-1,i,j+1,k)=m(n-1,i,j-1,k)+m(n-1,i,j+1,k)=m(n,i,j,k).$$

Thus for $k=4$ the matrices look like

$$\begin{pmatrix} a(0) & a(1) & a(2) & a(3) & a(4) \\ a(1) & a(0)+a(2) & a(1)+a(3) & a(2)+a(4) & a(3) \\ a(2) & a(1)+a(3) & a(0)+a(2)+a(4) & a(1)+a(3) & a(2) \\ a(3) & a(2)+a(4) & a(1)+a(3) & a(0)+a(2) & a(1) \\ a(4) & a(3) & a(2) & a(1) & a(0) \end{pmatrix}$$

Now we can prove that $A_{n,k}$ and $B_{n,k}$ have the same size.

**Corollary 1.3**

*For all* $n\in\mathbb{N}$ *and* $k\ge 1$ *we have* $\left|A_{n,k}\right|=\left|B_{n,k}\right|$.

Proof

We know already that

$$\left|A_{2n,k}\right| = m\left(2n, \left\lfloor \frac{k+1}{2} \right\rfloor, \left\lfloor \frac{k+1}{2} \right\rfloor, k\right)$$

$$\left|A_{2n+1,k}\right| = m\left(2n+1, \left\lfloor \frac{k+1}{2} \right\rfloor, \left\lfloor \frac{k-1}{2} \right\rfloor, k\right).$$

By (1.10) we get

$$m(n,i,j,k) = \sum_{\ell=0}^{i} m(n,0,j-i+2\ell,k)$$

$$m\left(2n, \left\lfloor \frac{k+1}{2} \right\rfloor, \left\lfloor \frac{k+1}{2} \right\rfloor, k\right) = \sum_{\ell=0}^{\left\lfloor \frac{k+1}{2} \right\rfloor} m(2n,0,2\ell,k) = \sum_{j=0}^{k} m(2n,0,j,k) = \left|B_{2n,k}\right|,$$

$$m\left(2n+1, \left\lfloor \frac{k+1}{2} \right\rfloor, \left\lfloor \frac{k-1}{2} \right\rfloor, k\right) = \sum_{\ell=0}^{\left\lfloor \frac{k+1}{2} \right\rfloor} m(2n+1,1,2\ell,k) = \sum_{j=0}^{k} m(2n+1,0,j,k) = \left|B_{2n+1,k}\right|.$$

The next result gives in combination with Corollary 1.3 another proof of (1.1).

**Proposition 1.4**

$$\left|B_{n,k}\right| = \sum_{\ell\in\mathbb{Z}} (-1)^{\ell} \binom{n}{\left\lfloor \frac{n+(k+2)\ell}{2} \right\rfloor}.$$

The following proof uses an idea by S.V. Ault and Ch. Kicey [2].

**Proof**

We show that for $-1 \le j \le k+1$

$$m(n,0,j,k) = \sum_{\ell\in\mathbb{Z}} \binom{n}{\left\lfloor \frac{j+n+1}{2} \right\rfloor + (k+2)\ell} - \sum_{\ell\in\mathbb{Z}} \binom{n}{\left\lfloor \frac{j+n+2}{2} \right\rfloor + (k+2)\ell}. \qquad (1.11)$$

To show this formula it suffices to check the recursion and the initial and boundary values.

$$m(0,0,j,k) = \sum_{\ell\in\mathbb{Z}} \binom{0}{\left\lfloor \frac{j+1}{2} \right\rfloor + (k+2)\ell} - \sum_{\ell\in\mathbb{Z}} \binom{0}{\left\lfloor \frac{j+2}{2} \right\rfloor + (k+2)\ell} = [j=0].$$

Since

$$\binom{n}{\left\lfloor \frac{n+j+1}{2} \right\rfloor + (k+2)\ell} = \binom{n}{n - \left\lfloor \frac{n+j+1}{2} \right\rfloor - (k+2)\ell} = \binom{n}{\left\lfloor \frac{n-j}{2} \right\rfloor - (k+2)\ell}$$

we get

$$m(n,0,-1,k) = \sum_{j\in\mathbb{Z}} \binom{n}{\left\lfloor \frac{n}{2} \right\rfloor + (k+2)j} - \sum_{j\in\mathbb{Z}} \binom{n}{\left\lfloor \frac{n+1}{2} \right\rfloor + (k+2)j} = 0$$

and

$$m(n,0,k+1,k) = \sum_{\ell\in\mathbb{Z}} \binom{n}{\left\lfloor \frac{k+2+n}{2} \right\rfloor + (k+2)\ell} - \sum_{\ell\in\mathbb{Z}} \binom{n}{\left\lfloor \frac{k+2+n+1}{2} \right\rfloor + (k+2)\ell}$$
$$= \sum_{\ell\in\mathbb{Z}} \binom{n}{\left\lfloor \frac{k+2+n+1}{2} \right\rfloor - (k+2)(\ell+1)} - \sum_{\ell\in\mathbb{Z}} \binom{n}{\left\lfloor \frac{k+2+n}{2} \right\rfloor - (k+2)(\ell+1)} = -m(n,0,k+1,k)$$

and thus $m(n,0,k+1,k) = 0.$

The recurrence follows from

$$\binom{n}{\left\lfloor \frac{j+n}{2} \right\rfloor + (k+2)\ell} = \binom{n-1}{\left\lfloor \frac{j-1+n-1}{2} \right\rfloor + (k+2)\ell} + \binom{n-1}{\left\lfloor \frac{j+1+n-1}{2} \right\rfloor + (k+2)\ell}.$$

Therefore

$$\sum_{j=-1}^{k} m(n,0,j,k) = \sum_{\ell\in\mathbb{Z}} \binom{n}{\left\lfloor \frac{n}{2} \right\rfloor + (k+2)\ell} - \sum_{\ell\in\mathbb{Z}} \binom{n}{\left\lfloor \frac{n+k+2}{2} \right\rfloor + (k+2)\ell}$$
$$= \sum_{\ell\in\mathbb{Z}} (-1)^{\ell} \binom{n}{\left\lfloor \frac{n+(k+2)\ell}{2} \right\rfloor}.$$

**Remark 1.3**

Proposition 1.4 can also be deduced from general results about lattice paths in corridors. E.g. [16], formula (9) implies that the number of walks on $P_{k+1}$ from $1$ to $m$ is $0$ if $n - m \equiv 0 \bmod 2$ and else $\sum_{j\in\mathbb{Z}} \binom{n}{\frac{m+n-1}{2} + (k+2)j} - \sum_{j\in\mathbb{Z}} \binom{n}{\frac{m+n+1}{2} + (k+2)j}.$

Both results can be combined to give (1.11).

### 1.3. Some useful facts about the matrices $M_k$

Let us give some more information about the matrices $M_k$.

To this end and for later applications we recall some facts about Fibonacci and Lucas polynomials.

The Fibonacci polynomials $F_n(x,s)=\sum_{k=0}^{\left\lfloor \frac{n-1}{2}\right\rfloor}\binom{n-1-k}{k}x^{n-1-2k}s^k$ satisfy the recurrence relation $F_n(x,s)=xF_{n-1}(x,s)+sF_{n-2}(x,s)$ with initial values $F_0(x,s)=0$ and $F_1(x,s)=1$ and the Lucas polynomials $L_n(x,s)=\sum_{k=0}^{\left\lfloor \frac{n}{2}\right\rfloor}\binom{n-k}{k}\frac{n}{n-k}s^k x^{n-2k}$ satisfy $L_n(x,s)=xL_{n-1}(x,s)+sL_{n-2}(x,s)$ with initial values $L_0(x,s)=2$ and $L_1(x,s)=x.$

Most identities about these polynomials can easily be proved by using the well-known Binet formulae

$F_n(x,s)=\frac{\alpha^n-\beta^n}{\alpha-\beta}$ and $L_n(x,s)=\alpha^n+\beta^n$ if $\alpha=\alpha(x,s)=\frac{x+\sqrt{x^2+4s}}{2}$ and $\beta=\beta(x,s)=\frac{x-\sqrt{x^2+4s}}{2}$ are the roots of the equation $z^2-xz-s=0.$

Let us do this for some formulae which will be needed in the sequel:

The identity

$$L_n(x,s)=F_{n+1}(x,s)+sF_{n-1}(x,s) \tag{1.12}$$

follows from $\left(\alpha^n+\beta^n\right)\left(\alpha-\beta\right)=\alpha^{n+1}-\beta^{n+1}-\alpha\beta\left(\alpha^{n-1}-\beta^{n-1}\right),$ the identity

$$F_{2n}(x,s)=F_n(x,s)L_n(x,s) \tag{1.13}$$

from $\alpha^{2n}-\beta^{2n}=\left(\alpha^n-\beta^n\right)\left(\alpha^n+\beta^n\right)$ and

$$F_{k+1}(x,s)^2+sF_k(x,s)^2=F_{2k+1}(x,s) \tag{1.14}$$

from $\left(\alpha^{k+1}-\beta^{k+1}\right)^2-\alpha\beta\left(\alpha^k-\beta^k\right)^2=\left(\alpha-\beta\right)\left(\alpha^{2k+1}-\beta^{2k+1}\right).$

Since $\alpha\left(x+y,-xy\right)=x$ and $\beta\left(x+y,-xy\right)=y$ we get the well-known identities

$$L_n\left(x+y,-xy\right)=x^n+y^n,$$
$$F_n(x+y,-xy)=\frac{x^n-y^n}{x-y}. \tag{1.15}$$

If we choose $x=e^{\frac{j\pi}{k+2}i},\ y=e^{-\frac{j\pi}{k+2}i}$ for $1\le j\le k+1$ we get

$F_{k+2}\left(2\cos\frac{j\pi}{k+2},-1\right)=0$ or since $F_{k+2}(x,-1)$ is a monic polynomial of degree $k$

$$F_{k+2}\left(x,-1\right)=\prod_{j=1}^{k+1}\left(x-2\cos\frac{j\pi}{k+2}\right). \tag{1.16}$$

The Fibonacci polynomials can be represented as the determinant

$$F_{k+2}\left(x,-1\right)=\det\begin{pmatrix} x & -1 & 0 & 0 & \cdots & 0\\ -1 & x & -1 & 0 & \cdots & 0\\ 0 & -1 & x & -1 & \cdots & 0\\ 0 & 0 & -1 & x & \cdots & 0\\ \vdots & \vdots & \vdots & \vdots & \ddots & \vdots\\ 0 & 0 & 0 & 0 & \cdots & x \end{pmatrix} \tag{1.17}$$

which follows immediately from their recurrence relation.

The right-hand side can be interpreted as the characteristic polynomial of the matrix $M_k$.

This is one of the reasons why Fibonacci polynomials play such a dominant role in this field.

I became aware of this fact through the blog post [27] by Qiaochu Yuan.

By (1.16) the eigenvalues of $M_k$ are given by $\lambda_j=2\cos\frac{j\pi}{k+2}$ for $1\le j\le k+1$.

Then $v_j=\left(F_1(\lambda_j,-1),F_2(\lambda_j,-1),\cdots,F_{k+1}(\lambda_j,-1)\right)^t$ is an eigenvector corresponding to $\lambda_j$.

For $M_k v_j=\lambda_j v_j$ is equivalent with $F_{\ell-1}\left(\lambda_j,-1\right)+F_{\ell+1}\left(\lambda_j,-1\right)=\lambda_j F_\ell\left(\lambda_j,-1\right)$ for $1\le\ell\le k+1$.

Note that $F_0\left(\lambda_j,-1\right)=F_{k+2}\left(\lambda_j,-1\right)=0$.

Since by (1.15) $F_\ell(\lambda_j,-1)=F_\ell(2\cos\frac{j\pi}{k+2},-1)=\dfrac{e^{\frac{\ell j\pi i}{k+2}}-e^{-\frac{\ell j\pi i}{k+2}}}{e^{\frac{j\pi i}{k+2}}-e^{-\frac{j\pi i}{k+2}}}=\dfrac{\sin\frac{\ell j\pi}{k+2}}{\sin\frac{j\pi}{k+2}}$

the eigenvectors are (up to scaling) given by $v_j = \left( \sin\frac{j\pi}{k+2}, \sin\frac{2j\pi}{k+2}, \cdots, , \sin\frac{(k+1)j\pi}{k+2} \right)^t$.

The normalized eigenvectors are $\sqrt{\frac{2}{k+2}} v_j$ since

$$\sum_{\ell=1}^{k+1} \left( \sin\frac{\ell j\pi}{k+2} \right)^2 = -\frac{1}{4}\sum_{\ell=1}^{k+1} \left( e^{\frac{2\ell j\pi}{k+2}i} - 2 + e^{-\frac{2\ell j\pi}{k+2}i} \right) = -\frac{1}{4}\left(-4-2k\right) = \frac{k+2}{2}.$$

Since $M_k$ is obviously symmetric we see that the matrix

$$U = \left( \sqrt{\frac{2}{k+2}} v_1, \sqrt{\frac{2}{k+2}} v_2, \cdots, \sqrt{\frac{2}{k+2}} v_{k+1} \right)$$

is orthogonal. Let $\Lambda_k = \left( \lambda_j [i=j] \right)_{i,j=1}^{k+1}$ be the diagonal matrix whose entries are the eigenvalues. Then $M_k = U\Lambda_k U^{-1} = U\Lambda_k U^t$.

Therefore from $M_k^n = U\Lambda_k^n U^t$ we get the known trigonometric representation

$$m(n,0,j,k) = \frac{2}{k+2}\sum_{\ell=1}^{k+1} \sin\frac{\ell\pi}{k+2} \sin\frac{\ell j\pi}{k+2} \left( 2\cos\frac{j\pi}{k+2} \right)^n. \tag{1.18}$$

References may be found in the recent paper [14] by Stefan Felsner and Daniel Heldt where similar results are obtained and in the survey article [17] by Christian Krattenthaler.

### 1.4. Generating functions

The generating functions of these number sequences turn out to be quotients of Fibonacci and Lucas polynomials or equivalently quotients of Chebyshev polynomials.

In the same way as above we see that $\det\left(I_k - M_k x\right) = F_{k+2}(1,-x^2)$.

From $\left(I_k - M_k x\right)^{-1} = \sum_{n\ge 0} M_k^n x^n$ and Cramer's Rule $\left(I_k - M_k x\right)^{-1} = \frac{adj(I_k - M_k x)}{\det\left(I_k - M_k x\right)}$ we find by considering the top-left entry of these matrices that the generating function $v_k(x)$ of the numbers $m(n,0,0,k)$ is given by

$$v_k(x) = \sum_{n\ge 0} m(n,0,0,k)x^n = \sum_{n\ge 0} m(2n,0,0,k)x^n = \frac{F_{k+1}\left(1,-x^2\right)}{F_{k+2}\left(1,-x^2\right)}. \tag{1.19}$$

Helmut Prodinger has kindly brought my attention to the paper [3] by N.G. de Bruijn, D.E. Knuth and S.O. Rice which gives another approach to this formula.

The numbers $m(2n,0,0,k)$ count the Dyck paths with height $\le k$.

Such a path $P$ is either the trivial path $(0,0) \to (0,0)$ of length $0$ or has a uniquely determined decomposition $P = UP_1DUP_2D\cdots UP_jD$ where each $P_i$ is a Dyck path with height $\le k-1$.

Therefore $v_k(x)$ satisfies

$$v_k(x) = 1 + x^2 v_{k-1}(x) + \left(x^2 v_{k-1}(x)\right)^2 + \left(x^2 v_{k-1}(x)\right)^3 + \cdots = \frac{1}{1 - x^2 v_{k-1}(x)}. \tag{1.20}$$

Note that $v_0(x) = 1$.

For arbitrary Dyck paths (1.20) gives the well-known fact that

$$v(x) = v_\infty(x) = \frac{1}{1 - x^2 v(x)} \quad \text{or} \quad v(x) = \frac{1 - \sqrt{1-4x^2}}{2x^2} = \sum_{n \ge 0} C_n x^{2n} \quad \text{where } C_n = \frac{1}{n+1}\binom{2n}{n}$$

is a Catalan number.

This implies that $m(2n,0,0,k) = C_n$ for $n \le k$.

From (1.20) we deduce the well-known generating function of bounded Dyck paths

$$v_k(x) = \frac{F_{k+1}(1,-x^2)}{F_{k+2}(1,-x^2)}. \tag{1.21}$$

For this holds for $k = 0$. If it is true for $k-1$ then

$$v_k(x) = \frac{1}{1 - x^2 v_{k-1}(x)} = \frac{1}{1 - \dfrac{x^2 F_k(1,-x^2)}{F_{k+1}(1,-x^2)}}$$
$$= \frac{F_{k+1}(1,-x^2)}{F_{k+1}(1,-x^2) - x^2 F_k(1,-x^2)} = \frac{F_{k+1}(1,-x^2)}{F_{k+2}(1,-x^2)}.$$

Let now more generally

$$v_k(x,j) = \sum_{n \ge 0} m(n,0,j,k) x^n. \tag{1.22}$$

Since each path $P$ from $(0,0)$ to $(n,j)$ has a unique decomposition

$P = P_0 U P_1 U P_2 \cdots U P_j$ where $P_\ell$ is a Dyck path bounded by $k - \ell$ we get

$$v_k(x,j) = x^j v_k(x) v_{k-1}(x) \cdots v_{k-j}(x). \tag{1.23}$$

By (1.23) and (1.21) we get

$$v_k(x,j) = x^j \frac{F_{k-j+1}\left(1,-x^2\right)}{F_{k+2}\left(1,-x^2\right)}. \tag{1.24}$$

As shown in [6] and [10] the generating functions of the sequences $\left(a(n,k)\right)_{n\geq 0}$ are given by

$$\sum_{n\geq 0} a(n,2k+1)x^n = \frac{F_{k+1}(1,-x^2)}{F_{k+2}(1,-x^2) - xF_{k+1}(1,-x^2)} \tag{1.25}$$

and

$$\sum_{n\geq 0} a(n,2k)x^n = \frac{F_{k+1}(1,-x^2) + xF_k(1,-x^2)}{L_{k+1}(1,-x^2)}. \tag{1.26}$$

Observe that $\deg\left(F_{k+2}(1,-x^2) - xF_{k+1}(1,-x^2)\right) = k+1$ and $\deg L_{k+1}(1,-x^2) = 2\left\lfloor \frac{k+1}{2} \right\rfloor$.

Since by (1.14) and (1.13)

$$\frac{F_{k+1}(1,-x^2)}{F_{k+2}(1,-x^2) - xF_{k+1}(1,-x^2)} = \frac{F_{k+1}(1,-x^2)\left(F_{k+2}(1,-x^2) + xF_{k+1}(1,-x^2)\right)}{F_{2k+3}(1,-x^2)}$$

and

$$\frac{F_{k+1}(1,-x^2) + xF_k(1,-x^2)}{L_{k+1}(1,-x^2)} = \frac{F_{k+1}(1,-x^2)\left(F_{k+1}(1,-x^2) + xF_k(1,-x^2)\right)}{F_{2k+2}\left(1,-x^2\right)}$$

both formulae (1.25) and (1.26) can be written compactly as

$$\sum_{n\geq 0} a(n,k)x^n = \frac{F_{\left\lfloor \frac{k+2}{2} \right\rfloor}\left(1,-x^2\right)\left(F_{\left\lfloor \frac{k+3}{2} \right\rfloor}\left(1,-x^2\right) + xF_{\left\lfloor \frac{k+1}{2} \right\rfloor}\left(1,-x^2\right)\right)}{F_{k+2}\left(1,-x^2\right)}. \tag{1.27}$$

Let us give another proof of these formulae.

As above $\left(I_k - R_k x\right)^{-1} = \sum_{n\geq 0} R_k^n x^n$ and Cramer's Rule gives $\left(I_k - R_k x\right)^{-1} = \dfrac{adj(I_k - R_k x)}{\det\left(I_k - R_k x\right)}$.

$$\det\left(I_k - xR_k\right) = \det\begin{pmatrix} 1-x & -x & 0 & \cdots & 0 & 0 \\ -x & 1 & -x & \cdots & 0 & 0 \\ 0 & -x & 1 & \cdots & 0 & 0 \\ \vdots & \vdots & \vdots & \ddots & \vdots & \vdots \\ 0 & 0 & 0 & \cdots & 1 & -x \\ 0 & 0 & 0 & \cdots & -x & 1 \end{pmatrix} = F_{k+2}\left(1,-x^2\right) - xF_{k+1}\left(1,-x^2\right).$$

For by expanding with respect to the first column we get
$(1-x)F_{k+1}\left(1,-x^2\right) - x^2F_k\left(1,-x^2\right) = F_{k+2}\left(1,-x^2\right) - xF_{k+1}\left(1,-x^2\right).$

Thus we get again

$$\sum_{n\geq 0} a(n,2k+1)x^n = \frac{F_{k+1}\left(1,-x^2\right)}{F_{k+2}\left(1,-x^2\right) - xF_{k+1}\left(1,-x^2\right)}.$$

Let $\alpha = \alpha\left(1,-x^2\right)$ and $\beta = \beta\left(1,-x^2\right)$. Then the right-hand side can also be written as

$$\frac{\alpha^{k+1} - \beta^{k+1}}{\alpha^{k+2} - \beta^{k+2} - x\left(\alpha^{k+1} - \beta^{k+1}\right)}.$$

By Corollary 1.2 we get $\left(I_k - S_k x\right)^{-1} = \sum_{n\geq 0} S_k^n x^n$ and $\left(I_k - S_k x\right)^{-1} = \dfrac{adj(I_k - S_k x)}{\det\left(I_k - S_k x\right)}.$

$$\det\left(I_k - S_k x\right) = \det\begin{pmatrix} 1 & -x & 0 & \cdots & 0 & 0 \\ -2x & 1 & -x & \cdots & 0 & 0 \\ 0 & -x & 1 & \cdots & 0 & 0 \\ \vdots & \vdots & \vdots & \ddots & \vdots & \vdots \\ 0 & 0 & 0 & \cdots & 1 & -x \\ 0 & 0 & 0 & \cdots & -x & 1 \end{pmatrix} = L_{k+1}\left(1,-x^2\right),$$

because $F_{k+1}\left(1,-x^2\right) - 2x^2F_k\left(1,-x^2\right) = F_{k+2}\left(1,-x^2\right) - x^2F_k\left(1,-x^2\right) = L_{k+1}\left(1,-x^2\right).$

By considering the upper-left-most entry we get

$$\sum_{n\geq 0} a(2n,2k)x^{2n} = \frac{F_{k+1}\left(1,-x^2\right)}{L_{k+1}\left(1,-x^2\right)} = \frac{1}{\alpha-\beta}\frac{\alpha^{k+1} - \beta^{k+1}}{\alpha^{k+1} + \beta^{k+1}}.$$

For $k \to \infty$ this converges to the well-known formula $\sum_{n\geq 0}\binom{2n}{n}x^n = \dfrac{1}{\alpha-\beta} = \dfrac{1}{\sqrt{1-4x^2}}.$

By considering the entry $(1,1)$ we get in the same way

$$\sum_{n\geq 0} a(2n+1,2k)x^{2n} = \frac{F_k\left(1,-x^2\right)}{L_{k+1}\left(1,-x^2\right)}.$$

Thus (1.26) is proved.

Another derivation of (1.27) is due to Helmut Prodinger (personal communication):

Since $a_k(x)=\sum_{j=0}^{k} v_k(x,j)$ we have to show that

$$\frac{F_{\left\lfloor\frac{k+2}{2}\right\rfloor}\left(1,-x^2\right)\left(F_{\left\lfloor\frac{k+3}{2}\right\rfloor}\left(1,-x^2\right)+xF_{\left\lfloor\frac{k+1}{2}\right\rfloor}\left(1,-x^2\right)\right)}{F_{k+2}\left(1,-x^2\right)} = \sum_{j=0}^{k} x^j \frac{F_{k-j+1}\left(1,-x^2\right)}{F_{k+2}\left(1,-x^2\right)}. \qquad (1.28)$$

This is equivalent with

$$F_{k+1}^2\left(1,-x^2\right) = \sum_j x^{2j} F_{2k-2j+1}\left(1,-x^2\right),$$

$$F_{k+1}\left(1,-x^2\right)F_k\left(1,-x^2\right) = \sum_j x^{2j} F_{2k-2j}\left(1,-x^2\right).$$

The first identity reduces to

$$\sum_{j=0}^{k}(\alpha\beta)^j \frac{\alpha^{2k-2j+1}-\beta^{2k-2j+1}}{\alpha-\beta} = \frac{1}{\alpha-\beta}\left(\alpha^{2k+1}\sum_{j=0}^{k}\left(\frac{\beta}{\alpha}\right)^j - \beta^{2k+1}\sum_{j=0}^{k}\left(\frac{\alpha}{\beta}\right)^j\right)$$

$$= \frac{\alpha\left(\alpha^{2k+1}-\beta^{k+1}\alpha^k\right)+\beta\left(\beta^{2k+1}-\alpha^{k+1}\beta^k\right)}{(\alpha-\beta)^2} = \left(\frac{\alpha^{k+1}-\beta^{k+1}}{\alpha-\beta}\right)^2$$

and the second one to

$$\sum_{j=0}^{k}(\alpha\beta)^j \frac{\alpha^{2k-2j}-\beta^{2k-2j}}{\alpha-\beta} = \frac{1}{\alpha-\beta}\left(\alpha^{2k}\sum_{j=0}^{k}\left(\frac{\beta}{\alpha}\right)^j - \beta^{2k}\sum_{j=0}^{k}\left(\frac{\alpha}{\beta}\right)^j\right)$$

$$= \frac{\alpha\left(\alpha^{2k}-\beta^{k+1}\alpha^{k-1}\right)+\beta\left(\beta^{2k}-\alpha^{k+1}\beta^{k-1}\right)}{(\alpha-\beta)^2} = \left(\frac{\alpha^{k+1}-\beta^{k+1}}{\alpha-\beta}\right)\left(\frac{\alpha^k-\beta^k}{\alpha-\beta}\right).$$

## 2. Polynomials associated with $A_{n,k}$.

### 2.1. Definitions and known results

Instead of the numbers $\left|A_{n,k}\right|$ we consider the following weights. Define a peak as a vertex preceded by an up-step $U$ and followed by a down-step $D$, and a valley as a vertex preceded by a down- step $D$ and followed by an up-step $U$. The height of a vertex is its $y-$coordinate. The peaks with a height at least $1$ and the valleys with height at most $-2$ are called extremal points. Let $E(v)$ be the set of $x-$coordinates of the extremal points of the path $v, e(v)=\left|E(v)\right|$ the number of extremal points of $v$ and $\iota(v)=\sum_{i\in E(v)} i$ the sum of the $x-$coordinates of the extremal points.

Following [20] we defined in [7] and [8] the weight of $v$ by $w(v,q)=q^{\iota(v)}t^{e(v)}$ and considered the polynomials $w\left(A_{n,k},q\right)=\sum_{v\in A_{n,k}} q^{\iota(v)}t^{e(v)}$. These polynomials are intimately connected with Rogers-Ramanujan type theorems.

In the present paper we consider only the case $q=1$ and study the polynomials

$$a(n,k,t)=\sum_{v\in A_{n,k}} t^{e(v)} \tag{2.1}$$

in more detail. It is obvious that $\deg\left(a(n,k,t)\right)=\left\lfloor\frac{n}{2}\right\rfloor$ for $k>1$ because the maximal degree is obtained by the path $UDUD\cdots$.

If we set $\binom{n}{k}=0$ for $n<0$ it follows from the results in [7] and [8] that for $k\geq 1$ these polynomials can be written in the following form:

$$a(n,k,t)=\sum_{j\in\mathbb{Z}}(-1)^j\sum_{\ell\geq|j|}\binom{\left\lfloor\frac{n+(k-2)j}{2}\right\rfloor}{\ell-j}\binom{\left\lfloor\frac{n+1-(k-2)j}{2}\right\rfloor}{\ell+j}t^\ell \tag{2.2}$$

For $t=1$ we have of course $a(n,k,1)=\left|A_{n,k}\right|$.

**Remark 2.1**

A direct proof that (2.2) implies

$$a(n,k,1)=\sum_{j\in\mathbb{Z}}(-1)^j\binom{n}{\left\lfloor\frac{n+(k+2)j}{2}\right\rfloor} \tag{2.3}$$

follows from the fact that $\left\lfloor \frac{n+kj}{2} \right\rfloor + \left\lfloor \frac{n+1-kj}{2} \right\rfloor = n.$

For

$$\sum_{\ell=|j|}^{\infty} \binom{\left\lfloor \frac{n+(k-2)j}{2} \right\rfloor}{\ell-j} \binom{\left\lfloor \frac{n+1-(k-2)j}{2} \right\rfloor}{\ell+j} = \sum_{\ell=|j|}^{\infty} \binom{\left\lfloor \frac{n+(k-2)j}{2} \right\rfloor}{\ell-j} \binom{n-\left\lfloor \frac{n+(k-2)j}{2} \right\rfloor}{n-\left\lfloor \frac{n+(k-2)j}{2} \right\rfloor-\ell-j}$$

$$= \sum_{i=-\infty}^{\infty} \binom{\left\lfloor \frac{n+(k-2)j}{2} \right\rfloor}{i} \binom{n-\left\lfloor \frac{n+(k-2)j}{2} \right\rfloor}{n-\left\lfloor \frac{n+(k-2)j}{2} \right\rfloor-i-2j}$$

$$= \sum_{i=-\infty}^{\infty} \binom{\left\lfloor \frac{n+(k-2)j}{2} \right\rfloor}{i} \binom{n-\left\lfloor \frac{n+(k-2)j}{2} \right\rfloor}{\left\lfloor \frac{n+1-(k+2)j}{2} \right\rfloor-i} = \binom{n}{\left\lfloor \frac{n+1-(k+2)j}{2} \right\rfloor} = \binom{n}{\left\lfloor \frac{n+(k+2)j}{2} \right\rfloor}.$$

For $n \le k$ we have $a(n,k,t) = \sum_{\ell \ge 0} \binom{\left\lfloor \frac{n}{2} \right\rfloor}{\ell} \binom{\left\lfloor \frac{n+1}{2} \right\rfloor}{\ell} t^{\ell}.$

The simplest special cases are

$a(n,1,t) = 1,$

$a(n,2,t) = (1+t)^{\left\lfloor \frac{n}{2} \right\rfloor}.$

As a generalization of (1.5) we get

$$a(n,3,t) = F_{n+1}(1,t) = \sum_{k=0}^{\left\lfloor \frac{n}{2} \right\rfloor} \binom{n-k}{k} t^k = \sum_{j\in\mathbb{Z}} (-1)^j \sum_{\ell \ge |j|} \binom{\left\lfloor \frac{n+j}{2} \right\rfloor}{\ell-j} \binom{\left\lfloor \frac{n+1-j}{2} \right\rfloor}{\ell+j} t^{\ell}. \qquad (2.4)$$

The first terms of $a(n,3,t)$ are

$\left(1,1,1+t,1+2t,1+3t+t^2,1+4t+3t^2,1+5t+6t^2+t^3,1+6t+10t^2+4t^3,\cdots\right)$

To show that $a(n,3,t) = F_{n+1}(1,t)$ consider a path in $A_{n,1}$. If the next to the last point is not extremal then the path is the unique continuation of a path in $A_{n-1,1}$, if it is extremal then the last two steps are a peak or a valley and the rest of the path belongs to $A_{n-2,1}$. Therefore we have $a(n,1,t) = a(n-1,1,t) + ta(n-2,1,t)$. The initial values are $a(0,1,t) = a(1,1,t) = 1$.

## 2.2. Generating functions for the polynomials $a(n,k,t)$.

### 2.2.1.

We now want to determine the generating functions for the polynomials $a(n,k,t)$.

To this end we introduce the polynomials

$$\Phi_n(x,t) = F_n\left(1+(1-t)x^2,-x^2\right) - x^2 F_{n-1}\left(1+(1-t)x^2,-x^2\right). \quad (2.5)$$

They satisfy $\Phi_n(x,t) = \left(1+(1-t)x^2\right)\Phi_{n-1}(x,t) - x^2\Phi_{n-2}(x,t)$ with initial values $\Phi_0(x,t) = \Phi_1(x,t) = 1.$ For $t=1$ we get $\Phi_n(x,1) = F_{n+1}\left(1,-x^2\right).$

**Theorem 2.1**

*For $k \geq 0$ the generating function of the sequence $\left(a(n,2k+1,t)\right)_{n\geq 0}$ is*

$$\sum_{n\geq 0} a(n,2k+1,t)x^n = \frac{\Phi_k(x,t)}{\Phi_{k+1}(x,t) - x\Phi_k(x,t)}. \quad (2.6)$$

For the proof we first observe that for the polynomials $a(n,2k+1,t)$ Proposition 1.1 remains true if we define the weight of a non-negative lattice path $p$ by $w(p) = t^{e(p)}$ where $e(p)$ denotes the number of peaks of $p$.

**Proposition 2.1**

*The polynomial $a(n,2k+1,t)$ is the weight of all non-negative lattice paths starting from $(0,0)$ and ending in $(n,0)$, where besides up-steps and down-steps also horizontal moves $(i,0) \to (i+1,0)$ on height $0$ are allowed and the maximal height of a path is $k$.*

**Proof of Theorem 2.1**

Let $A(k,x,t)$ be the generating function of these lattice paths and let $D(k,x,t)$ be the generating function of all Dyck paths with maximal height $k$ with the same weight.

Then

$$A(k,x,t) = 1 + xA(k,x,t) + x^2 tA(k,x,t) + x^2\left(D(k-1,x,t)-1\right)A(k,x,t)$$

because such a path is either trivial or begins with a horizontal step or with $UD$ or with $UPD$ where $P$ is a non-trivial Dyckpath of height $\leq k-1$.

Therefore we get

$$A(k,x,t) = \frac{1}{1-x-x^2t+x^2-x^2D(k-1,x,t)}.$$

Since $D(0,x,t)=1$ we get $A(1,x,t)=\frac{1}{1-x-x^2t}.$

In order to compute $A(k,x,t)$ we must first compute $D(k,x,t).$

Here we have

$D(k,x,t)=1+x^2tD(k,x,t)+x^2\left(D(k-1,x,t)-1\right)$ and thus

$$D(k,x,t)=\frac{1}{1-x^2t+x^2-x^2D(k-1,x,t)}. \tag{2.7}$$

This gives

$$D(k,x,t)=\frac{\Phi_k(x,t)}{\Phi_{k+1}(x,t)}. \tag{2.8}$$

This follows by induction because

$$\frac{1}{1+(1-t)x^2-x^2\frac{\Phi_{k-1}(x,t)}{\Phi_k(x,t)}}=\frac{\Phi_k(x,t)}{\left(1+(1-t)x^2\right)\Phi_k(x,t)-x^2\Phi_{k-1}(x,t)}=\frac{\Phi_k(x,t)}{\Phi_{k+1}(x,t)}.$$

Then we get

$$A(k,x,t)=\frac{1}{1+(1-t)x^2-x-x^2\frac{\Phi_{k-1}(x,t)}{\Phi_k(x,t)}}=\frac{\Phi_k(x,t)}{\Phi_{k+1}(x,t)-x\Phi_k(x,t)}.$$

If we let $k\to\infty$ in (2.7) we get the generating function $D(x,t)$ of all Dyck paths with weight given by the number of peaks $w(p)=t^{e(p)}$ (cf. [13])

$$D(x,t)=\frac{1+(1-t)x^2-\sqrt{\left(1+(1-t)x^2\right)^2-4x^2}}{2x^2}=\sum_{n\geq 0}N_n(t)x^{2n}. \tag{2.9}$$

Here $N_n(t)=\frac{1}{n}\sum_{k=0}^{n}\binom{n}{k}\binom{n}{k-1}t^k$ for $n>0$ and $N_0(t)=1$ denotes a Narayana polynomial (cf. OEIS A001263 or [18]). Note that there is a certain ambiguity of notation in the literature. In Petersen [21] the polynomials $C_n(t)=\frac{N_n(t)}{t}$ for $n>0$ with $C_0(t)=1$ are called Narayana polynomials. In [18] these polynomials are called associated Narayana polynomials.

If we let $k \to \infty$ in $A(k,x,t)$ we get the generating function

$$\sum_{n\geq 0}\left(\sum_{\ell\geq 0}\binom{\left\lfloor\frac{n}{2}\right\rfloor}{\ell}\binom{\left\lfloor\frac{n+1}{2}\right\rfloor}{\ell}t^{\ell}\right)x^n=\frac{1}{1+(1-t)x^2-x-x^2D(x,t)}. \tag{2.10}$$

In an analogous manner we introduce the polynomials

$$\Lambda_n(x,t)=L_n(1+(1-t)x^2,-x^2)-x^2L_{n-1}(1+(1-t)x^2,-x^2). \tag{2.11}$$

They satisfy $\Lambda_n(x,t)=\left(1+(1-t)x^2\right)\Lambda_{n-1}(x,t)-x^2\Lambda_{n-2}(x,t)$ with initial values $\Lambda_0(x,t)=1-(1-t)x^2$ and $\Lambda_1(x,t)=1-(1+t)x^2$. For $t=1$ we get $\Lambda_n(x,1)=L_{n+1}\left(1,-x^2\right)$.

**Theorem 2.2**

*For $k\geq 1$ the generating function of $\left(a(n,2k,t)\right)_{n\geq 0}$ is*

$$\sum_{n\geq 0}a(n,2k,t)x^n=\frac{\Phi_k(x,t)+x\Phi_{k-1}(x,t)}{\Lambda_k(x,t)}. \tag{2.12}$$

**Proof**

Let $B(k,x,t)=\sum_{n\geq 0}a(2n,2k,t)x^{2n}$.

Then $B(k,x,t)=1+x^2tB(k,x,t)+x^2B(k,x,t)+2x^2\left(D(k-1,x,t)-1\right)B(k,x,t)$

because each path is either empty or begins with $UD$ or with $DU$ or has the form $UPD$ or $DP'U$ where $P$ is a Dyck path with height $\leq k-1$ and $P'$ a reflected Dyck path with height $\leq k-1$.

Therefore

$$B(k,x,t)=\frac{1}{1-x^2t+x^2-2x^2D(k-1,x,t)}=\frac{1}{1+(1-t)x^2-2x^2\frac{\Phi_{k-1}(x,t)}{\Phi_k(x,t)}}$$
$$=\frac{\Phi_k(x,t)}{\Phi_{k+1}(x,t)-x^2\Phi_{k-1}(x,t)}=\frac{\Phi_k(x,t)}{\Lambda_k(x,t)}.$$

Let $C(k,x,t)=\sum_n a(2n+1,2k,t)x^{2n+1}$.

Then $C(k,x,t)=xD(k-1,x,t)B(k,x,t)$ because each path ends with $DP'$ where $P'$ is a reflected path of $P\in D(k-1,x,t)$. Therefore we get

$$C(k,x,t)=x\frac{\Phi_{k-1}(x,t)}{\Phi_k(x,t)}\frac{\Phi_k(x,t)}{\Lambda_k(x,t)}=x\frac{\Phi_{k-1}(x,t)}{\Lambda_k(x,t)}$$

If we let $k\to\infty$ in the formula

$$\sum_{n\geq 0}a(2n,2k,t)x^{2n}=\frac{\Phi_k(x,t)}{\Lambda_k(x,t)}=\frac{1}{1-x^2t+x^2-2x^2D(k-1,x,t)}$$ and observe (2.9) we get

$$\sum_{n\geq 0}\left(\sum_{k=0}^{n}\binom{n}{k}^2 t^k\right)x^{2n}=\frac{1}{\sqrt{\left(1+(1-t)x^2\right)^2-4x^2}}. \tag{2.13}$$

**2.2.2.**

Let us now consider some special cases.

The simplest special cases are

$$\sum_{n\geq 0}a(n,2,t)x^n=\frac{1+x}{(1+(1-t)x^2)-2x^2}=\frac{1+x}{1-(1+t)x^2}=\sum_{n\geq 0}(1+t)^{\left\lfloor\frac{n}{2}\right\rfloor}x^n$$

and

$$\sum_{n\geq 0}a(n,3,t)x^n=\frac{1}{(1+(1-t)x^2)-x(x+1)}=\frac{1}{1-x-tx^2}=\sum_{n\geq 0}F_{n+1}(1,t)x^n.$$

Let us also consider $a(n,4,t)$.

The first terms are

$1,1,1+t,1+2t,1+4t+t^2,1+5t+3t^2,1+7t+9t^2+t^3,1+8t+14t^2+4t^3,\cdots$.

Here we have

$$\sum_{n\geq 0}a(n,4,t)x^n=\frac{1+x-tx^2}{1-x^2-2tx^2-tx^4+t^2x^4}=\frac{1-tx^2}{1-x^2-2tx^2-tx^4+t^2x^4}+x\frac{1}{1-x^2-2tx^2-tx^4+t^2x^4}$$
$$=\sum_{n\geq 0}a(2n,4,t)x^{2n}+x\sum_{n\geq 0}a(2n+1,4,t)x^{2n}.$$

Since $\dfrac{1-tx^2}{1-x^2-2tx^2-tx^4+t^2x^4}-1=\dfrac{x^2+tx^2(1-tx^2)+tx^4}{1-x^2-2tx^2-tx^4+t^2x^4}$

we get by comparing coefficients and setting $a(-1,4,t)=0$

$a(2n+3,4,t) = a(2n+2,4,t) + ta(2n+1,4,t)$
$a(2n+2,4,t) = a(2n+1,4,t) + ta(2n,4,t) + ta(2n-1,4,t)$

Thus the polynomials $a(n,4,t)$ can easily be computed.

Let now $b(n,t) = a(2n,4,t^2) + ta(2n-1,4,t^2) = \sum_{k=0}^{2n} b_{n,k} t^k$. Then the table $\left(b_{n,k}\right)$ begins as follows

$$\begin{pmatrix} &&&&&1&&&&& \\ &&&&1&1&1&&&& \\ &&&1&1&4&2&1&&& \\ &&1&1&7&5&9&3&1&& \\ &1&1&10&8&26&14&16&4&1& \\ 1&1&13&11&52&34&70&30&25&5&1 \end{pmatrix} \tag{2.14}$$

The sum of the rows is $b(n,1) = a(2n,4,1) + a(2n-1,4,1) = 2\cdot 3^{n-1} + 3^{n-1} = 3^n$ and the alternating sum is

$$\sum_{k=0}^{2n} (-1)^k b_{n,k} = a(2n,4,1) - a(2n-1,4,1) = 2\cdot 3^{n-1} - 3^{n-1} = 3^{n-1}.$$

The table $\left(b_{n,k}\right)$ satisfies

$$\begin{aligned} b_{n,2k+1} &= b_{n-1,2k-1} + b_{n-1,2k}, \\ b_{n,2k} &= b_{n-1,2k-2} + 2b_{n-1,2k-1} + b_{n-1,2k}. \end{aligned} \tag{2.15}$$

It seems that $b_{2n,2n} = \sum_{k=0}^{n} \binom{n}{k}\binom{2n+k}{k}$ and $b_{2n+1,2n+1} = \sum_{k=0}^{n} \binom{n}{k}\binom{2n+k+1}{k}$.

### 2.3. Generating functions of the coefficients

**Theorem 2.3**

*For each* $k \ge 0$ *there exist polynomials* $u_j(x, 2k+1)$ *with integer coefficients such that*

$$\sum_{n\ge 0} a(n,2k+1,t)x^n = \frac{1}{1-x} + \frac{x^2 u_1(x,2k+1)}{(1-x)^2}t + \frac{x^4 u_2(x,2k+1)}{(1-x)^3}t^2 + \frac{x^6 u_3(x,2k+1)}{(1-x)^4}t^3 + \cdots. \tag{2.16}$$

*Especially we have* $u_1(x,2k+1) = \sum_{j=1}^{k} x^{2j-2}$ *and* $u_j(1,2k+1) = k^j$.

**Proof**

By induction it is easily verified that

$$\Phi_n(x,t) = 1 + tx^2 c_1\left(x^2\right) + t^2x^4 c_2\left(x^2\right) + \cdots \text{ with } c_1\left(x^2\right) = -\sum_{j=1}^{n}(n-j)x^{2j-2} \text{ and}$$

$\deg c_j(x^2) = 2n-2.$

Therefore $\Phi_{k+1}(x,t) - x\Phi_k(x,t) = 1 - x + tx^2 p_1(x) + t^2x^4 p_2(x) + \cdots + t^k x^{2k} p_k(x)$ for some polynomials $p_j(x)$ with $\deg p_j = 2k.$

By multiplying both sides of (2.16) with $\Phi_{k+1}(x,t) - x\Phi_k(x,t)$ we get

$$\left(\Phi_{k+1}(x,t) - x\Phi_k(x,t)\right)\sum_{n\geq 0} a(n,2k+1,t)x^n$$
$$= \left(1 - x + tx^2 p_1(x) + t^2x^4 p_2(x) + \cdots + t^k x^{2k} p_k(x)\right)$$
$$*\left(\frac{1}{1-x} + \frac{x^2 u_1(x,2k+1)}{(1-x)^2} t + \frac{x^4 u_2(x,2k+1)}{(1-x)^3} t^2 + \frac{x^6 u_3(x,2k+1)}{(1-x)^4} t^3 + \cdots\right) = \Phi_k(x,t).$$

Since the highest power of $t$ in $\Phi_k(x,t)$ is $t^{k-1}$ we see that the coefficient of $t^n$ vanishes for $n \geq k.$

Comparing coefficients of $t^n$ we get

$$\frac{u_n(x,2k+1)}{(1-x)^n} + \frac{u_{n-1}(x,2k+1)}{(1-x)^n} p_1(x) + \cdots + \frac{u_{n-k}(x,2k+1)}{(1-x)^{n+1-k}} p_k(x) = c_n\left(x^2\right)$$

where $c_n\left(x^2\right) = 0$ for $n \geq k.$

This is equivalent with

$$u_n(x,k) = -p_1(x)u_{n-1}(x,2k+1) - (1-x)p_2(x)u_{n-2}(x,2k+1) - \cdots + c_n\left(x^2\right)(1-x)^n \quad (2.17)$$

Therefore $u_n(x,k)$ is a polynomial with integer coefficients.

Let us compute $u_1(x,2k+1) = -p_1(x) + (1-x)c_1\left(x^2\right).$

We get $p_1(x) = -\sum_{j=1}^{n}(k-j)x^{2j-2} + x\sum_{j=1}^{n}(k-1-j)x^{2j-2}.$

This implies

$$u_1(x,2k+1) = -p_1(x) + (1-x)c_1\left(x^2\right) = \sum_{j=1}^{k}(k+1-j)x^{2j-2} - x\sum_{j=1}^{k}(k-j)x^{2j-2} - (1-x)\sum_{j=1}^{k}(k-j)x^{2j-2}$$
$$= \sum_{j=1}^{k}(k+1-j)x^{2j-2} - \sum_{j=1}^{k}(k-j)x^{2j-2} = \sum_{j=1}^{k}x^{2j-2}.$$

Thus $u_1(1,2k+1)=k$ and by (2.17) we get $u_j(1,2k+1)=k^j$.

**Some special cases.**

For $k=0$ we get $u_j(x,0)=0.$

For $k=1$ we have $\Phi_2(x,t)-x\Phi_1(x,t)=1-x-tx^2$ and therefore $p_1(x)=-1.$

This implies $u_n(x,1)=1$ for all $n.$ Thus we get

$$\sum_{n\geq 0}a(n,3,t)x^n=\frac{1}{1-x-tx^2}=\sum_{j\geq 0}t^j\frac{x^{2j}}{(1-x)^{j+1}}.$$

Let us also consider $k=2.$ Here we get $\Phi_2(x,t)=1-tx^2$ and

$$\Phi_3(x,t)-x\Phi_2(x,t)=1-x+tx^2\left(-2+x-x^2\right)+t^2x^4$$

This implies $u_1(x,5)=\left(2-x+x^2\right)-(1-x)=1+x^2$ and

$$u_j(x,5)=(2-x+x^2)u_{j-1}(x,5)-(1-x)u_{j-2}(x,5) \qquad (2.18)$$

for $j\geq 2.$

**Theorem 2.4**

*For $k\geq 1$*

$$\sum_{n\geq 0}a(n,2k,t)x^n=\frac{1}{1-x}+\frac{x^2v_1(x,2k)}{(1-x)^2(1+x)}t+\frac{x^4v_2(x,2k)}{(1-x)^3(1+x)^2}t^2+\frac{x^6v_3(x,2k)}{(1-x)^4(1+x)^3}t^3+\cdots$$

*for some polynomials $v_j(x,2k)$ with integer coefficients. Moreover $v_1(x,2k)=\sum_{j=0}^{2k-2}x^j$ and*

$$\begin{aligned}v_j(1,2k)&=(2k-1)^j,\\ v_j(-1,2k)&=(2k-1)^{j-1}.\end{aligned} \qquad (2.19)$$

**Proof**

By multiplying both sides with $\Lambda_k(x,t)=1-x^2+tx^2h_1\left(x^2\right)+\cdots+t^kx^{2k}h_k\left(x^2\right)$ we get

$$\begin{aligned}&\Lambda_k(x,t)\sum_{n\geq 0}a(n,2k,t)x^n\\ &=\left(1-x^2+tx^2h_{k,1}\left(x^2\right)+\cdots+t^kx^{2k}h_{k,k}\left(x^2\right)\right)\left(\frac{1}{1-x}+\frac{x^2v_1(x,2k)}{(1-x)^2(1+x)}t+\frac{x^4v_2(x,2k)}{(1-x)^3(1+x)^2}t^2+\frac{x^6v_3(x,2k)}{(1-x)^4(1+x)^3}t^3+\cdots\right)\\ &=\Phi_k(x,t)+x\Phi_{k-1}(x,t).\end{aligned}$$

As above we see that each $v_j(x,2k)$ is a polynomial with integer coefficients and that $v_1(x,2k)=\sum_{j=0}^{2k-2} x^j$. Since $\deg h_{k,j}=2k-2j$ it follows that $\deg v_j=(2k-2)j$.

By induction we see that $h_{k,1}\left(x^2\right)=-k-x^2-x^4-\cdots-x^{2k-2}$. Therefore we get that the leading coefficient of $v_j(x,2k)$ is 1 and that $v_j(\pm1,2k)=h_{k,1}(1)v_{j-1}(\pm1,2k)$ which gives $v_j(1,2k)=(2k-1)^j$ and $v_j(-1,2k)=(2k-1)^{j-1}$.

**Remark**

If we set $v_j(x,2k+1)=(1+x)^j u_j(x,k)$ then we have for all $k\ge1$

$$\sum_{n\ge0} a(n,k,t)x^n = \frac{1}{1-x}+\frac{x^2v_1(x,k)}{(1-x)^2(1+x)}t+\frac{x^4v_2(x,k)}{(1-x)^3(1+x)^2}t^2+\frac{x^6v_3(x,k)}{(1-x)^4(1+x)^3}t^3+\cdots \tag{2.20}$$

Then $v_j(x,k)$ is a polynomial with degree $(k-2)j$ and satisfies $v_j(1,k)=(k-1)^j$.

Note that for each $k$ we have $v_1(x,k)=\sum_{j=0}^{k-2} x^j$ and that the leading coefficient of $v_j(x,k)$ is 1. It seems that moreover each $v_j(x,k)$ has non-negative coefficients.

**Some special cases**

For $k=3$ we have $v_j(x,3)=(1+x)^j$.

The next case is more interesting.

We know already that

$$\sum_{n\ge0} a(n,4,t)x^n = \frac{1+x-tx^2}{1-x^2-2tx^2-tx^4+t^2x^4}.$$

We already know that

$$\sum_{n\ge0} a(n,4,t)x^n = \frac{1}{1-x}+\frac{x^2v_1(x,4)}{(1-x)^2(1+x)}t+\frac{x^4v_2(x,4)}{(1-x)^3(1+x)^2}t^2+\frac{x^6v_3(x,4)}{(1-x)^4(1+x)^3}t^3+\cdots \tag{2.21}$$

for some polynomials $v_j(x,4)$.

Multiplying both sides by $1-x^2-2tx^2-tx^4+t^2x^4$ and comparing coefficients of $t^j$ we get

$$v_j(x,4)=(2+x^2)v_{j-1}(x,4)-(1-x^2)v_{j-2}(x,4). \tag{2.22}$$

The initial values are $v_0(x,4)=1$ and $v_1(x,4)=1+x+x^2$ by direct computation.

If we compute the polynomials $v_j(x,4)=\sum_{\ell=0}^{2j} c(j,\ell)x^\ell$ we get the following array of the coefficients $c(j,\ell)$

| | | | | | | | | | | |
|---|---|---|---|---|---|---|---|---|---|---|
| 1 | | | | | | | | | | |
| 1 | 1 | 1 | | | | | | | | |
| 1 | 2 | 4 | 1 | 1 | | | | | | |
| 1 | 3 | 9 | 5 | 7 | 1 | 1 | | | | |
| 1 | 4 | 16 | 14 | 26 | 8 | 10 | 1 | 1 | | |
| 1 | 5 | 25 | 30 | 70 | 34 | 52 | 11 | 13 | 1 | 1 |

(2.22) implies the following formulae:

$c(-1,\ell)=c(-2,\ell)=0, c(0,\ell)=[\ell=0],\ c(1,\ell)=[\ell\le 2],$

$c(j,\ell)=2c(j-1,\ell)-c(j-2,\ell)+c(j-1,\ell-2)+c(j-2,\ell-2).$

This implies that $c(n,2n)=c(n,2n-1)=1$ and $c(n,\ell)\ge c(n-1,\ell).$ Thus all coefficients are non-negative.

Surprisingly this is almost the same array as (2.14) . More precisely

$v_j(x,4)=x^{2j}b\left(j,\frac{1}{x}\right)$ since both sides satisfy the same recurrence (2.22).

Thus we have $v_j(1,4)=\sum_{\ell=0}^{2j} c(j,\ell)=3^j$ and $v_j(-1,4)=\sum_{\ell=0}^{2j}(-1)^\ell c(j,\ell)=3^{j-1}.$

By (2.18) $v_j(x,5)$ satisfies

$v_j(x,5)=\left(x^3+x+2\right)v_{j-1}(x,5)+(1-x)(1+x)^2 v_{j-2}(x,5)$

with $v_0(x,5)=1$ and $v_1(x,5)=(1+x)\left(1+x^2\right).$

For $v_j(x,6)$ we get

$v_j(x,6)=\left(x^4+x^2+3\right)v_{j-1}(x,6)+\left(2x^4+x^2-3\right)v_{j-2}(x,6)+\left(1-x^2\right)^2 v_{j-3}(x,6).$

If we let $k \to \infty$ we get

$$a(n,t)=\sum_{\ell\geq 0}\binom{\left\lfloor \frac{n}{2}\right\rfloor}{\ell}\binom{\left\lfloor \frac{n+1}{2}\right\rfloor}{\ell}t^{\ell}. \tag{2.23}$$

Here we have

$$\sum_{n\geq 0}a(n,t)x^{n}=\frac{1}{1-x}+\frac{x^{2}r_{0}(x)}{(1-x)^{3}(1+x)}t+\frac{x^{4}r_{1}(x)}{(1-x)^{5}(1+x)^{3}}t^{2}+\frac{x^{6}r_{2}(x)}{(1-x)^{7}(1+x)^{5}}t^{3}+\cdots \tag{2.24}$$

where

$$r_{j}(x)=\sum_{\ell=0}^{j}\binom{j}{\ell}^{2}x^{2\ell}+\sum_{\ell=1}^{j}jN_{j,\ell}x^{2\ell-1}=\sum_{\ell=0}^{j}\binom{j}{\ell}^{2}x^{2\ell}+\sum_{\ell=1}^{j}\binom{j}{\ell}\binom{j}{\ell-1}x^{2\ell-1}. \tag{2.25}$$

The numbers $N_{n,k}=\binom{n}{k-1}\binom{n}{k}\frac{1}{n}=\binom{n}{k-1}\binom{n-1}{k-1}\frac{1}{k}$ are the Narayana numbers (cf. OEIS A001263) and the polynomials $N_{n}(x)=\sum_{k=0}^{n}N_{n,k}x^{k}$ with $N_{0}(x)=1$ the Narayana polynomials.

The coefficient table of $\left(r_{j}(x)\right)_{j\geq 0}$ is ( OEIS A247644)

$$\begin{array}{ccccccccc} & & & & 1 & & & & \\ & & & 1 & 1 & 1 & & & \\ & & 1 & 2 & 4 & 2 & 1 & & \\ & 1 & 3 & 9 & 9 & 9 & 3 & 1 & \\ 1 & 4 & 16 & 24 & 36 & 24 & 16 & 4 & 1 \end{array} \tag{2.26}$$

Note that (2.24) is equivalent with

$$(1-x)^{2}\left(1-x^{2}\right)^{2k-1}\sum_{n\geq 0}\binom{\left\lfloor \frac{n+2k}{2}\right\rfloor}{k}\binom{\left\lfloor \frac{n+1+2k}{2}\right\rfloor}{k}x^{n}=r_{k-1}(x) \tag{2.27}$$

for all $k$.

To show this identity we make use of the identities

$$(1-x)^{2k+1}\sum_{n\geq 1}\binom{n+k-1}{k}\binom{n+k}{k}x^{n-1}=\sum_{j=1}^{k}\binom{k-1}{j-1}\binom{k+1}{j}x^{j-1} \tag{2.28}$$

and

$$(1-x)^{2k+1}\sum_{n\geq 0}\binom{n+k}{k}^2 x^n=\sum_{j=0}^{k}\binom{k}{j}^2 x^j. \tag{2.29}$$

Identity (2.28) can also be formulated as

$$(1-x)^{2k+1}\sum_{n\geq 0}N_{n+k,k+1}x^n=N_k(x). \tag{2.30}$$

These imply

$$\left(1-x^2\right)^{2k+1}\sum_{n\geq 0}\binom{\left\lfloor\frac{n+2k}{2}\right\rfloor}{k}\binom{\left\lfloor\frac{n+1+2k}{2}\right\rfloor}{k}x^n$$
$$=\left(1-x^2\right)^{2k+1}\sum_{n\geq 1}\binom{n+k-1}{k}\binom{n+k}{k}x^{2n-1}+\left(1-x^2\right)^{2k+1}\sum_{k\geq 0}\binom{n+k}{k}^2 x^{2n}$$
$$=\sum_{j=0}^{k}\binom{k}{j}^2 x^{2j}+\sum_{j=1}^{k}\binom{k-1}{j-1}\binom{k+1}{j}x^{2j-1}=(1+x)^2 r_{k-1}(x)$$

by observing that

$$(1+x)^2\left(\sum_{j=0}^{k-1}\binom{k-1}{j}^2 x^{2j}+\sum_{j=1}^{k-1}\binom{k-1}{j-1}\binom{k-1}{j}x^{2j-1}\right)=\sum_j\binom{k-1}{j}^2 x^{2j}+2\sum_j\binom{k-1}{j-1}^2 x^{2j-1}+\sum_j\binom{k-1}{j-1}^2 x^{2j}$$
$$+\sum_j\binom{k-1}{j-1}\binom{k-1}{j}x^{2j-1}+2\sum_j\binom{k-1}{j-1}\binom{k-1}{j}x^{2j}+\sum_j\binom{k-1}{j-2}\binom{k-1}{j-1}x^{2j-1}$$
$$=\sum_j\left(\binom{k-1}{j-1}+\binom{k-1}{j}\right)^2 x^{2j}+\sum_j\binom{k-1}{j-1}\left(2\binom{k-1}{j-1}+\binom{k-1}{j}+\binom{k-1}{j-2}\right)$$
$$=\sum_j\binom{k}{j}^2 x^{2j}+\sum_j\binom{k-1}{j-1}\binom{k+1}{j}x^{2j-1}.$$

To prove (2.28) and (2.29) we show more generally that for $m \geq 0$

$$\left(1-x\right)^{2k+1} \sum_{n\geq 0} \binom{n+k}{k} \binom{n+k-m}{k} x^{n-m} = \sum_{j=0}^{k} \binom{k-m}{j} \binom{2k-j}{k} (-1)^j (1-x)^j$$
$$= \sum_{j=m}^{k} \binom{k-m}{j-m} \binom{k+m}{j} x^{j-m}.$$

Since $\sum_{n\geq 0} \binom{n+k}{k} x^n = \frac{1}{(1-x)^{k+1}}$

we get

$$\sum_{n\geq 0} \binom{n+k}{k} \binom{n+k-m}{k} x^{n-m} = \frac{1}{k!} D^k \frac{x^{k-m}}{(1-x)^{k+1}} = \frac{1}{k!} D^k \frac{\left(1-(1-x)\right)^{k-m}}{(1-x)^{k+1}} = \frac{1}{k!} D^k \sum_{j=0}^{k-m} (-1)^j \binom{k-m}{j} \left(1-x\right)^{j-}$$
$$= \sum_{j=0}^{k-m} (-1)^j \binom{k-m}{j} \binom{2k-j}{k} \left(1-x\right)^{j-2k-1}$$

It remains to show that

$$\sum_{j=0}^{k-m} (-1)^j \binom{k-m}{j} \binom{2k-j}{k} \left(1-x\right)^j = \sum_{j=m}^{k} \binom{k-m}{j-m} \binom{k+m}{j} x^{j-m}. \tag{2.31}$$

This follows by comparing coefficients in two different expansions of $(1+z)^{k+m}(x+z)^{k-m}$.

On the one hand we have

$$(1+z)^{k+m}(x+z)^{k-m} = (1+z)^{k+m}(x-1+1+z)^{k-m} = \sum_{j} \binom{k-m}{j} (x-1)^j (1+z)^{2k-j}$$
$$= \sum_{j,\ell} \binom{k-m}{j} (x-1)^j \binom{2k-j}{\ell} z^\ell$$

On the other hand we get

$$(1+z)^{k+m}(x+z)^{k-m} = \sum_{j,\ell} \binom{k+m}{j} \binom{k-m}{\ell} x^\ell z^{j+k-m-\ell}.$$

Comparing the coefficients of $z^k$ in both sums we get (2.31).

**Remark 2.2**

For $m = 0$ identity (2.31) reduces to

$$\sum_{j=0}^{k} (-1)^j \binom{k}{j} \binom{2k-j}{k} \left(1-x\right)^j = \sum_{j=0}^{k} \binom{k}{j}^2 x^j. \tag{2.32}$$

This identity is mentioned in OEIS A063007 without proof. A combinatorial proof has been given by H.S. Wilf [26], p. 117. The above proof is inspired by the paper [24] by Jocelyn Quaintance, which contains tables of seven unpublished manuscript notebooks of H. W. Gould from 1945 – 1990. Similar identities can be found in [4] and the literature cited there.

A more general identity of this sort has been proved combinatorially in [4]. It can equivalently be formulated as

$$\sum_{j=0}^{n}\binom{n}{j}\binom{n+2m+x}{j+m}z^{j}=\sum_{j=0}^{n}\binom{n}{j}\binom{2n+2m+x-j}{n+m}(z-1)^{j} \tag{2.33}$$

and reduces to (2.31) for $(n,x)\to(k-m,0)$.

Since there seems to be some interest in such identities I will give another proof of (2.27).

We first show by induction that

$$b(n,j,x):=(1-x)^{k+j+1}\frac{D^{j}}{j!}\frac{x^{n}}{(1-x)^{k+1}}=\sum_{i=0}^{n}\binom{j+k-n}{k-i}\binom{n}{i}x^{i}. \tag{2.34}$$

Observe that $b(n,0,x)=(1-x)^{k+1}\dfrac{x^{n}}{(1-x)^{k+1}}=x^{n}=\sum_{i=0}^{n}\binom{j+k-n}{k-i}\binom{n}{i}x^{i}$

because $\binom{k-n}{k-i}=[i=n]$ for $i\le n$. We also have

$$b(0,j,x)=(1-x)^{k+j+1}\frac{D^{j}}{j!}\frac{1}{(1-x)^{k+1}}=\binom{j+k}{k}.$$

Since $D^{j}xf(x)=xD^{j}f(x)+jD^{j-1}f(x)$ the sequence $\left(b(n,j,x)\right)_{n\ge 0,j\ge 0}$ satisfies

$b(n,j,x)=xb(n-1,j,x)+(1-x)b(n-1,j-1,x).$

Comparing coefficients this is equivalent with

$$\binom{j+k-n}{k-i}\binom{n}{i}=\binom{j+k-n+1}{k-i+1}\binom{n-1}{i-1}+\binom{j+k-n}{k-i}\binom{n-1}{i}-\binom{j+k-n}{k-i+1}\binom{n-1}{i-1}.$$

This is clear because the right-hand side is

$$\binom{j+k-n+1}{k-i+1}\binom{n-1}{i-1}-\binom{j+k-n}{k-i+1}\binom{n-1}{i-1}+\binom{j+k-n}{k-i}\binom{n-1}{i}$$
$$=\binom{j+k-n}{k-i}\binom{n-1}{i-1}+\binom{j+k-n}{k-i}\binom{n-1}{i}=\binom{j+k-n}{k-i}\binom{n}{i}.$$

Thus for $m \le k$ (2.34) gives

$$(1-x)^{k+j+1}\frac{D^j}{j!}\frac{x^{k-m}}{(1-x)^{k+1}}=\sum_{i=0}^{k-m}\binom{j+m}{k-i}\binom{k-m}{i}x^i. \tag{2.35}$$

If we choose $m=0$ and $j=k$ we get

$$\begin{aligned}(1-x)^{2k+1}\sum_{n\ge 0}\binom{n+k}{k}^2 x^n &= (1-x)^{2k+1}\frac{D^k}{k!}\sum_{n\ge 0}\binom{n+k}{k}x^{n+k}\\ =(1-x)^{2k+1}\frac{D^k}{k!}\frac{x^k}{(1-x)^{k+1}} &= \sum_{j=0}^{k}\binom{k}{j}^2 x^j.\end{aligned} \tag{2.36}$$

For $m=1$ and $j=k-2$ we get

$$\begin{aligned}(1-x)^{2k+1}\sum_{n\ge 0}\binom{n+k}{k}\binom{n+k-1}{k-2}x^{n+1} &= (1-x)^{2k+1}\frac{D^{k-2}}{(k-2)!}\sum_{n\ge 0}\binom{n+k}{k}x^{n+k-1}\\ =(1-x)^{2k+1}\frac{D^{k-2}}{(k-2)!}\frac{x^{k-1}}{(1-x)^{k+1}} &= \sum_{i=0}^{k-1}\binom{k-1}{i-1}\binom{k-1}{i}x^i.\end{aligned} \tag{2.37}$$

Comparing coefficients gives

$$(1-x)^2\sum_{n\ge 0}\binom{\left\lfloor\frac{n+2k}{2}\right\rfloor}{k}\binom{\left\lfloor\frac{n+1+2k}{2}\right\rfloor}{k}x^n = (1-x)^2\sum_{n\ge 0}\binom{n+k}{k}^2 x^{2n}+(1-x)^2\sum_{n\ge 0}\binom{n+k}{k}\binom{n+k+1}{k}x^{2n+1}$$
$$=\sum_{n\ge 0}\binom{n+k-1}{k-1}^2 x^{2n}+\sum_{n\ge 0}\binom{n+k}{k}\binom{n+k-1}{k-2}x^{2n+1}.$$

Combining these identities we finally get the desired result (2.27).

Computing $b(n,j,x)$ in another way we get

$$\begin{aligned}(1-x)^{k+j+1}\frac{D^j}{j!}\frac{x^{k-m}}{(1-x)^{k+1}} &= (1-x)^{k+j+1}\frac{D^j}{j!}\frac{\sum_{\ell=0}^{k-m}\binom{k-m}{\ell}(-1)^\ell(1-x)^\ell}{(1-x)^{k+1}}\\ =(1-x)^{k+j+1}\frac{D^j}{j!}\sum_{\ell=0}^{k-m}\binom{k-m}{\ell}(-1)^\ell(1-x)^{\ell-k-1} &= \sum_{\ell=0}^{k-m}(-1)^\ell\binom{k-m}{\ell}\binom{k+j-\ell}{j}(1-x)^\ell.\end{aligned} \tag{2.38}$$

Comparison of (2.35) and (2.38) gives

$$\sum_{i=0}^{k-m}\binom{j+m}{k-i}\binom{k-m}{i}x^i=\sum_{\ell=0}^{k-m}(-1)^\ell\binom{k-m}{\ell}\binom{k+j-\ell}{j}(1-x)^\ell. \tag{2.39}$$

For $m=0$ and $j=k$ this reduces again to (2.32).

For $m=1$ and $j=k-2$ we get

$$\sum_{i=0}^{k-1}\binom{k-1}{k-i}\binom{k-1}{i}x^i=\sum_{\ell=0}^{k-1}(-1)^\ell\binom{k-1}{\ell}\binom{2k-2-\ell}{k-2}(1-x)^\ell,$$

which by $k\to k+1$ can be written as

$$\begin{aligned}\sum_{i=0}^{k}N_{k,i}x^i&=\sum_{i=0}^{k}\binom{k}{i-1}\binom{k}{i}\frac{1}{k}x^i=\sum_{\ell=0}^{k}(-1)^\ell\frac{1}{k+1}\binom{k+1}{\ell}\binom{2k-\ell}{k}(1-x)^\ell\\&=\sum_{\ell=0}^{k}(-1)^\ell\frac{1}{n-\ell+1}\binom{2n-2\ell}{n-\ell}\binom{2k-\ell}{k}(1-x)^\ell.\end{aligned} \tag{2.40}$$

In this form it has been proved in [5], (2.2) and [18], (1.3).

**Remark 2.3**

A slight modification of the above proof gives the following $q-$analogue of (2.27):

$$\left(x;q^k\right)_2\left(qx^2;q\right)_{2k-1}\sum_{n\geq 0}\begin{bmatrix}\left\lfloor\frac{n+2k}{2}\right\rfloor\\k\end{bmatrix}\begin{bmatrix}\left\lfloor\frac{n+1+2k}{2}\right\rfloor\\k\end{bmatrix}x^n=r_{k-1}(x,q) \tag{2.41}$$

with

$$r_n(x,q)=\sum_{j=0}^{2n}q^{\left\lfloor\frac{(j+1)^2}{4}\right\rfloor}\begin{bmatrix}n\\\left\lfloor\frac{j}{2}\right\rfloor\end{bmatrix}\begin{bmatrix}n\\\left\lfloor\frac{j+1}{2}\right\rfloor\end{bmatrix}x^j=\sum_{j=0}^{n}q^{2\binom{j+1}{2}}\begin{bmatrix}n\\j\end{bmatrix}^2x^{2j}+\sum_{j=1}^{n}q^{j^2}\begin{bmatrix}n\\j\end{bmatrix}\begin{bmatrix}n\\j-1\end{bmatrix}x^{2j-1}. \tag{2.42}$$

Here $\begin{bmatrix}n\\k\end{bmatrix}$ is $q-$binomial coefficient and $\left(x;q\right)_n=\prod_{j=0}^{n-1}\left(1-q^jx\right)$.

Comparing coefficients gives

$$\left(x;q^k\right)_2 \sum_{n\geq 0} \begin{bmatrix} \left[\frac{n+2k}{2}\right] \\ k \end{bmatrix} \begin{bmatrix} \left[\frac{n+1+2k}{2}\right] \\ k \end{bmatrix} x^n = \sum_{n\geq 0} q^n \begin{bmatrix} n+k-1 \\ k-1 \end{bmatrix}^2 x^{2n} + \sum_{n\geq 0} q^{n+1} \begin{bmatrix} n+k \\ k \end{bmatrix} \begin{bmatrix} n+k-1 \\ k-2 \end{bmatrix} x^{2n+1}.$$

Let $D_q$ be the $q-$differential operator defined by $D_q f(x) = \dfrac{f(x)-f(qx)}{(1-q)x}$ and define polynomials $b(n,j,x,q)$ by $b(n,j,x,q) = \left(x;q\right)_{k+j+1} \dfrac{D_q^j}{[j]!} \dfrac{x^n}{\left(x;q\right)_{k+1}}$.

Then

$$b(n,j,x,q) = \sum_{i=0}^{n} q^{i(j+i-n)} \begin{bmatrix} j+k-n \\ k-i \end{bmatrix} \begin{bmatrix} n \\ i \end{bmatrix} x^i. \tag{2.43}$$

Since $D_q^j x f(x) = q^j x D_q^j f(x) + [j] D_q^{j-1} f(x)$ the sequence $b(n,j,x,q)$ satisfies

$$b(n,j,x,q) = q^j x b(n-1,j,x,q) + (1-q^{k+j}x) b(n-1,j-1,x,q). \tag{2.44}$$

Comparing coefficients this is equivalent with

$$\begin{bmatrix} j+k-n \\ k-i \end{bmatrix} \begin{bmatrix} n \\ i \end{bmatrix} = q^{n-i} \begin{bmatrix} j+k-n+1 \\ k-i+1 \end{bmatrix} \begin{bmatrix} n-1 \\ i-1 \end{bmatrix} + \begin{bmatrix} j+k-n \\ k-i \end{bmatrix} \begin{bmatrix} n-1 \\ i \end{bmatrix} - q^{k-2i+n+1} \begin{bmatrix} j+k-n \\ k-i+1 \end{bmatrix} \begin{bmatrix} n-1 \\ i-1 \end{bmatrix}$$

.

The right-hand side is

$$q^{n-i} \begin{bmatrix} j+k-n+1 \\ k-i+1 \end{bmatrix} \begin{bmatrix} n-1 \\ i-1 \end{bmatrix} - q^{k-i+1} \begin{bmatrix} j+k-n \\ k-i+1 \end{bmatrix} \begin{bmatrix} n-1 \\ i-1 \end{bmatrix} + \begin{bmatrix} j+k-n \\ k-i \end{bmatrix} \begin{bmatrix} n-1 \\ i \end{bmatrix}$$

$$= q^{n-i} \begin{bmatrix} n-1 \\ i-1 \end{bmatrix} \begin{bmatrix} j+k-n \\ k-i \end{bmatrix} + \begin{bmatrix} j+k-n \\ k-i \end{bmatrix} \begin{bmatrix} n-1 \\ i \end{bmatrix} = \begin{bmatrix} n \\ i \end{bmatrix} \begin{bmatrix} j+k-n \\ k-i \end{bmatrix}.$$

Now observe that

$$\frac{1}{\left(x;q\right)_{k+1}} = \sum_{n\geq 0} \begin{bmatrix} n+k \\ k \end{bmatrix} x^n.$$

Therefore

$$\left(x;q\right)_{2k+1} \sum_{n\geq 0} \begin{bmatrix} n+k \\ k \end{bmatrix}^2 x^n = \left(x;q\right)_{2k+1} \frac{D^k}{[k]!} \frac{x^k}{\left(x;q\right)_{k+1}} = b(k,k,x,q) = \sum_{i=0}^{k} q^{i^2} \begin{bmatrix} k \\ i \end{bmatrix}^2 x^i$$

and

$$\left(qx^2;q\right)_{2k-1}\sum_{n\geq 0}\begin{bmatrix}n+k-1\\k-1\end{bmatrix}^2\left(qx^2\right)^n=\sum_{i=0}^{k-1}q^{i^2}\begin{bmatrix}k-1\\i\end{bmatrix}^2\left(qx^2\right)^i=\sum_{j=0}^{k-1}q^{j^2+j}\begin{bmatrix}k-1\\j\end{bmatrix}^2x^{2j}.$$

In the same way we get

$$\left(x;q\right)_{2k-1}\sum_{n\geq 0}\begin{bmatrix}n+k\\k\end{bmatrix}\begin{bmatrix}n+k-1\\k-2\end{bmatrix}x^{n+1}=\left(x;q\right)_{2k-1}\frac{D^{k-2}}{[k-2]!}\frac{x^{k-1}}{\left(x;q\right)_{k+1}}$$
$$=b(k-1,k-2,x,q)=\sum_{i=0}^{k-1}q^{i(i-1)}\begin{bmatrix}k-1\\k-i\end{bmatrix}\begin{bmatrix}k-1\\i\end{bmatrix}x^i$$

and therefore

$$\left(qx^2;q\right)_{2k-1}\sum_{n\geq 0}q^{n+1}\begin{bmatrix}n+k\\k\end{bmatrix}\begin{bmatrix}n+k-1\\k-2\end{bmatrix}x^{2n+1}=\sum_{i=1}^{k-1}q^{i^2}\begin{bmatrix}k-1\\k-i\end{bmatrix}\begin{bmatrix}k-1\\i\end{bmatrix}x^{2i-1}.$$

From the easily verified formulae

$$x^n=\sum_{\ell=0}^{n}(-1)^\ell\begin{bmatrix}n\\\ell\end{bmatrix}q^{\binom{\ell+1}{2}-n\ell}\left(x;q\right)_\ell \tag{2.45}$$

and

$$\frac{D_q^j}{[j]!}\frac{\left(x;q\right)_\ell}{\left(x;q\right)_{k+1}}=q^{\ell j}\begin{bmatrix}k+j-\ell\\j\end{bmatrix}\frac{\left(x;q\right)_\ell}{\left(x;q\right)_{k+j+1}} \tag{2.46}$$

we get

$$\frac{D_q^j}{[j]!}\frac{x^n}{\left(x;q\right)_{k+1}}=\frac{D_q^j}{[j]!}\frac{\sum_{\ell=0}^{n}(-1)^\ell\begin{bmatrix}n\\\ell\end{bmatrix}q^{\binom{\ell+1}{2}-n\ell}\left(x;q\right)_\ell}{\left(x;q\right)_{k+1}}=\frac{\sum_{\ell=0}^{n}(-1)^\ell q^{\binom{\ell+1}{2}+\ell(j-n)}\begin{bmatrix}n\\\ell\end{bmatrix}\begin{bmatrix}k+j-\ell\\j\end{bmatrix}\left(x;q\right)_\ell}{\left(x;q\right)_{k+j+1}}. \tag{2.47}$$

As special cases we get

$$\left(x;q\right)_{2k+1}\frac{D_q^k}{[k]!}\frac{x^k}{\left(x;q\right)_{k+1}}=\sum_{\ell=0}^{k}(-1)^j\begin{bmatrix}k\\\ell\end{bmatrix}\begin{bmatrix}2k-\ell\\k\end{bmatrix}q^{\binom{\ell+1}{2}}\left(x;q\right)_\ell$$

and

$$\left(x;q\right)_{2k-1}\frac{D^{k-2}}{[k-2]!}\frac{x^{k-1}}{\left(x;q\right)_{k+1}}=\sum_{\ell=0}^{k-1}(-1)^\ell\begin{bmatrix}k-1\\\ell\end{bmatrix}\begin{bmatrix}2k-2-\ell\\k-2\end{bmatrix}q^{\binom{\ell}{2}}\left(x;q\right)_\ell.$$

Comparing (2.47) with (2.43) we get as $q-$analogue of (2.39)

$$\sum_{i=0}^{n} q^{i(j+i-n)} \begin{bmatrix} j+k-n \\ k-i \end{bmatrix} \begin{bmatrix} n \\ i \end{bmatrix} x^i = \sum_{\ell=0}^{n} (-1)^\ell q^{\binom{\ell+1}{2}+\ell(j-n)} \begin{bmatrix} n \\ \ell \end{bmatrix} \begin{bmatrix} k+j-\ell \\ j \end{bmatrix} \left(x;q\right)_\ell . \tag{2.48}$$

The most interesting special cases are

$$\sum_{i=0}^{n} q^{i^2} \begin{bmatrix} n \\ i \end{bmatrix}^2 x^i = \sum_{\ell=0}^{n} (-1)^\ell q^{\binom{\ell+1}{2}} \begin{bmatrix} n \\ \ell \end{bmatrix} \begin{bmatrix} 2n-\ell \\ n \end{bmatrix} \left(x;q\right)_\ell \tag{2.49}$$

for $j = n = k$ and

$$\sum_{i=0}^{k-1} q^{i(i-1)} \begin{bmatrix} k-1 \\ k-i \end{bmatrix} \begin{bmatrix} k-1 \\ i \end{bmatrix} x^i = \sum_{\ell=0}^{k-1} (-1)^\ell q^{\binom{\ell+1}{2}-\ell} \begin{bmatrix} k-1 \\ \ell \end{bmatrix} \begin{bmatrix} 2k-2-\ell \\ k-2 \end{bmatrix} \left(x;q\right)_\ell$$

for $j = k-2$ and $n = k-1$. By substituting $k \to k+1$ this can be written as

$$\sum_{i=0}^{k} q^{i(i-1)} \frac{1}{[k]} \begin{bmatrix} k \\ i \end{bmatrix} \begin{bmatrix} k \\ i-1 \end{bmatrix} x^i = \sum_{\ell=0}^{k} (-1)^\ell q^{\binom{\ell}{2}} \frac{1}{[k+1]} \begin{bmatrix} k+1 \\ \ell \end{bmatrix} \begin{bmatrix} 2k-\ell \\ k \end{bmatrix} \left(x;q\right)_\ell , \tag{2.50}$$

where $\frac{1}{[k]} \begin{bmatrix} k \\ i \end{bmatrix} \begin{bmatrix} k \\ i-1 \end{bmatrix}$ is a $q-$analogue of a Narayana number.

**The polynomials $v_j(x,k)$.**

After this digression let us state some more properties of the polynomials $v_j(x,k)$.

By comparing (2.20) with (2.24) we see that $\lim_{k\to\infty} v_j(x,k) = \frac{r_{j-1}(x)}{(1-x)^j(1+x)^{j-1}}$.

Let us show more generally that $\left[x^\ell\right] v_j(x,k) = \left[x^\ell\right] \frac{r_{j-1}(x)}{(1-x)^j(1+x)^{j-1}}$ for $\ell \le k-2$ if we expand the right-hand side into a power series.

Let $a(n,t) = \lim_{k\to\infty} a(n,k,t)$. Then $\left[t^j\right] a(n,k,t) = \left[t^j\right] a(n,t)$ for $n \le k-2+2j$ since for these $n$ each path with $j$ extremal points belongs to $A_{n,k}$, because only if $n = k-2+2j$ there is a path with $j$ extremal points which reaches $y = \left\lfloor \frac{k}{2} \right\rfloor$, for example $\left(UD\right)^{j-1} U^{\left\lfloor \frac{k}{2} \right\rfloor} D^{\left\lfloor \frac{k}{2} \right\rfloor}$.

Now we get the desired result

$$\left[x^{\ell}\right]x^{2j}v_j(x,k)=\left[x^{\ell}\right](1-x)^{j+1}(1+x)^j\sum_{n\le 2j+k-2}\left[t^j\right]a(n,k,t)x^n$$

$$=\left[x^{\ell}\right](1-x)^{j+1}(1+x)^j x^{2j}\sum_{n\le 2j+k-2}\left[t^j\right]a(n,t)x^n=\left[x^{\ell}\right]x^{2j}\frac{r_{j-1}(x)}{(1-x)^j(1+x)^{j-1}}$$

for $\ell\le k-2+2j.$

For $j=1$ we get again that $v_1(x,k)=1+x+\cdots+x^{k-2}.$

We could also obtain this result by using formulae (2.6) and (2.12) and develop the right-hand sides into Taylor series. This can be done by expressing $\Phi_k$ and $\Lambda_k$ using $\alpha$ and $\beta.$

The coefficient table is

$$\begin{array}{ccccccc} & & & 1 & & & \\ & & 1 & & 1 & & \\ & 1 & & 1 & & 1 & \\ 1 & & 1 & & 1 & & 1 \end{array}$$

with generating function for the row sums

$$\frac{1}{(1-z)(1-xz)}=1+(1+x)z+(1+x+x^2)z^2+\cdots.$$

Since $\frac{x^2}{(1-x)^2(1+x)}=\sum_{n\ge 0}\left\lfloor\frac{n}{2}\right\rfloor x^n$ we see that

$$[t]a(n,k,t)=\sum_{j=0}^{k-2}\max\left(\left\lfloor\frac{n-j}{2}\right\rfloor,0\right)\text{ and }[t]a(n,t)=\sum_{j=0}^{n}\left\lfloor\frac{j}{2}\right\rfloor=\left\lfloor\frac{n}{2}\right\rfloor\left\lfloor\frac{n+1}{2}\right\rfloor.$$

Thus $[t]a(n,k,t)=\left\lfloor\frac{n}{2}\right\rfloor\left\lfloor\frac{n+1}{2}\right\rfloor$ for $n\le k-1$ and $[t]a(n,2k+1,t)=\sum_{j=0}^{2k-1}\left\lfloor\frac{n-j}{2}\right\rfloor=k(n-k)$

for $n\ge 2k+1$ and $[t]a(n,2k,t)=\sum_{j=0}^{2k-2}\left\lfloor\frac{n-j}{2}\right\rfloor=k(n-k)-\left\lfloor\frac{n-2k+1}{2}\right\rfloor$ for $n\ge 2k.$

For $j=2$ the Taylor series expansion gives

$$v_2(x,k)=\frac{1+x+x^2-x^{k-1}\left(k+2x+(2-k)x^2\right)+x^{2k-1}}{(1-x)^2(1+x)}.$$

This can also be written as

$$v_2(x,k+2)=1+\sum_{i=1}^{k}x^i(i+1+x+x^2+\cdots+x^i) \tag{2.51}$$

Let $v_2(x,k+2)=\sum_{\ell=0}^{2k} u_k(\ell)x^\ell$.

The coefficient table $\left(u_k(\ell)\right)_{k\geq 0}$ is

$$\begin{array}{ccccccccccc}
 & & & & & 1 & & & & & \\
 & & & & 1 & 2 & 1 & & & & \\
 & & & 1 & 2 & 4 & 1 & 1 & & & \\
 & & 1 & 2 & 4 & 5 & 2 & 1 & 1 & & \\
 & 1 & 2 & 4 & 5 & 7 & 2 & 2 & 1 & 1 & \\
1 & 2 & 4 & 5 & 7 & 8 & 3 & 2 & 2 & 1 & 1
\end{array} \tag{2.52}$$

Let

$$\frac{1+x+x^2}{(1-x)^2(1+x)}=1+2x+4x^2+5x^3+7x^4+8x^5+\cdots=\sum_{\ell\geq 0}c_\ell x^\ell.$$

Then we already know that $u_k(\ell)=c_\ell$ for $0\leq \ell\leq k$.

Let $c=\left(1,2,4,5,7,8,\cdots\right)$ and $d=\left(d_\ell\right)=\left(1,1,2,2,3,3,\cdots\right)$ with generating function $\frac{1}{(1+x)(1-x)^2}$. Then we have $u_k(2k-\ell)=d_\ell=\binom{\left\lfloor \frac{\ell+2}{2}\right\rfloor}{1}$ for $0\leq \ell\leq k-1$.

The sequence $c$ is the sequence of all positive integers which are not multiples of $3$. Thus $c_{2n}=3n+1$ and $c_{2n+1}=3n+2$.

The generating function of the rows is therefore

$$\sum_{k\geq 0} v_2(x,k+2)z^k=\frac{1-z^2x^3}{(1-z)(1-zx)^2(1-zx^2)}. \tag{2.53}$$

For $x=1$ this reduces to $\frac{1+z}{(1-z)^3}=\sum_{k\geq 0}(k+1)^2z^k$.

For $x=-1$ we get $\frac{1+z^2}{(1-z^2)^2}=\sum_{n\geq 0}(2k+1)z^{2k}$.

This is no surprise because by (2.19) we know already that $v_2(1,k+2)=(k+1)^2$ and $v_2(-1,2k+2)=2k+1$ and $v_2(-1,2k+1)=0$.

It is perhaps interesting that the coefficient table (2.52) has the form

$$\begin{array}{ccccccccccc} & & & & & 1 & & & & & \\ & & & & 1 & c_1 & 1 & & & & \\ & & & 1 & c_1 & c_2 & d_1 & 1 & & & \\ & & 1 & c_1 & c_2 & c_3 & d_2 & d_1 & 1 & & \\ & 1 & c_1 & c_2 & c_3 & c_4 & d_3 & d_2 & d_1 & 1 & \\ 1 & c_1 & c_2 & c_3 & c_4 & c_5 & d_4 & d_3 & d_2 & d_1 & 1 \end{array}$$

The sequences $c$ and $d$ are uniquely determined by $v_2(1,k+2)=(k+1)^2$, $v_2(-1,2k+2)=2k+1$ and $v_2(-1,2k+1)=0$.

Observe also that

$$c_0+\left(c_1+d_0\right)+(c_2+d_1)+\cdots+\left(c_n+d_{n-1}\right)=1+3+5+\cdots+(2n+1)=(n+1)^2,$$
$$c_0-\left(c_1+d_0\right)+(c_2+d_1)+\cdots+\left(c_{2n}+d_{2n-1}\right)=1-3+\cdots+(4n+1)=2n+1$$

and $\sum_{k=0}^{2n+1}(-1)^k c_k=-(n+1)$ and $\sum_{k=0}^{2n}(-1)^k d_k=n+1.$

In general it seems that

$$\left[x^{kj-\ell}\right]v_j\left(x,k+2\right)=\left[x^{\ell}\right]\frac{1}{(1-x)^{j+1}(1+x)^j}=\binom{\left\lfloor \frac{\ell+2j}{2}\right\rfloor}{j} \text{ for } 0\le \ell\le k-1.$$

For $j\geq 3$ the situation becomes more complicated.

With the help of a computer the Taylor expansion gives

$$\begin{aligned}&\left(1-x\right)^3\left(1+x\right)^2 v_3(x,2k+1)=1+2x+4x^2+2x^3+x^4\\&-x^{2k}\left((1+k)(1+2k)+4(1+k)x+(5-4k^2)x^2-4(k-1)x^3+(k-1)(2k-1)x^4\right)\\&+2x^{4k+1}\left(1+2k+x+(1-2k)x^2\right)-x^{6k+2}\end{aligned}$$

and

$$\begin{aligned}&\left(1-x\right)^3\left(1+x\right)^2 v_3(x,2k)=1+2x+4x^2+2x^3+x^4\\&-x^{2k-1}\left(k(2k+1)+2(1+2k)x+4\left(1+k-k^2\right)x^2+2(3-2k)x^3+(k-1)(2k-3)x^4\right)\\&+2x^{4k-1}\left(2k+x+2(k-1)x^2\right)-x^{6k-1}\end{aligned}$$

This can be simplified to give

$$(1-x)^3(1+x)^2 v_3(x,k)=1+2x+4x^2+2x^3+x^4$$
$$-x^{k-1}\left(\binom{k+1}{2}+2(1+k)x+\left(4+2k-k^2\right)x^2+(6-2k)x^3+\frac{6-5k+k^2}{2}x^4\right)$$
$$+2x^{2k-1}\left(k+x-(k-2)x^2\right)-x^{3k-1}.$$

The coefficient table of $v_3(x,k+2)=\left(u(k,j)\right)_{j=0}^{3k}$ starts with

| | | | | | | | | | | | | |
|---|---|---|---|---|---|---|---|---|---|---|---|---|
| 1 | | | | | | | | | | | | |
| 1 | 3 | 3 | 1 | | | | | | | | | |
| 1 | 3 | 9 | 5 | 7 | 1 | 1 | | | | | | |
| 1 | 3 | 9 | 15 | 12 | 10 | 9 | 3 | 1 | 1 | | | |
| 1 | 3 | 9 | 15 | 27 | 16 | 20 | 12 | 14 | 3 | 3 | 1 | 1 |

This table has some unusual properties.

Each column $[x^{k-j}]v_3(x,k+2)$ for $k\geq j$ is given by
$\left(c_n\right)_{n\geq 0}=\left(1,3,9,15,27,37,55,69,93,111,\cdots\right)$

with generating function

$$\sum_{n\geq 0} c_n x^n=\frac{r_2(x)}{(1-x)^3(1+x)^2}=\frac{1+2x+4x^2+2x^3+x^4}{(1-x)^3(1+x)^2}. \tag{2.54}$$

This implies that

$c_{2n}=1+3n+5n^2$ and $c_{2n-1}=5n^2-3n+1.$

On the right-hand side of the table each column $[x^{2k+j}]v_3(x,k+2)$ for $k\geq j$ equals

$\left(1,1,3,3,6,6,\cdots\right)=\left(\binom{\left\lfloor\frac{n+4}{2}\right\rfloor}{2}\right)_{n\geq 0}$. Its generating function is $\dfrac{1}{(1-x)^3(1+x)^2}$.

If we write the table in the form

$$\left(u(n,k)\right)_{n,k=0}^{\infty}=\begin{pmatrix} \color{green}{1} & & & & & & & & & & & & \\ 1 & \color{green}{3} & 3 & \color{red}{1} & & & & & & & & & \\ 1 & 3 & \color{green}{9} & 5 & 7 & \color{red}{1} & 1 & & & & & & \\ 1 & 3 & 9 & \color{green}{15} & 12 & 10 & 9 & \color{red}{3} & 1 & 1 & & & \\ 1 & 3 & 9 & 15 & \color{green}{27} & 16 & 20 & 12 & 14 & \color{red}{3} & 3 & 1 & 1 \end{pmatrix}$$

then the table between the green and the red column looks like

$$\begin{pmatrix} & & & 3 & & & \\ & & 5 & & 7 & & \\ & 12 & & 10 & & 9 & \\ 16 & & 20 & & 12 & & 14 \end{pmatrix} = \begin{pmatrix} & & & u(1,2) & & & \\ & & u(2,3) & & u(2,4) & & \\ & u(3,4) & & u(3,5) & & u(3,6) & \\ u(4,5) & & u(4,6) & & u(4,7) & & u(4,8) \end{pmatrix}.$$

It has some curious properties:

Let $R_k = \left(u\left(n+1, 2n+2-k\right) - u\left(n, 2n-k\right)\right)_{n \geq k+1}$ be the sequence of successive differences of the south-east columns $\left(u(n, 2n-k)\right)_{n \geq k+1}$.

Then

$$R_{2k} = \left(4k+4, 2k+2, 4k+5, 2k+2, 4k+6, 2k+2, \cdots\right)$$

and

$$R_{2k+1} = \left(4k+5, 2k+2, 4k+6, 2k+2, 4k+7, 2k+2, \cdots\right).$$

Let $L_k = \left(u\left(n+1, n+2+k\right) - u\left(n, n+k+1\right)\right)_{n \geq k+1}$ be the sequence of successive differences of the south-west columns $\left(u(n, n+k+1)\right)_{n \geq k+1}$.

Then

$$L_{2k} = \left(k+2, 5k+7, k+4, 5k+11, k+6, 5k+15, \cdots\right)$$

and

$$L_{2k-1} = \left(k+2, 5k+5, k+4, 5k+9, k+6, 5k+13, \cdots\right).$$

The exact values are

$$u(2k-1, 2k+2j) = 3k^2 - \frac{j(5j+1)}{2} \text{ for } 0 \leq j \leq k-1$$

and

$$u(2k-1, 2k-1+2j) = 3k^2 - \frac{j(5j-1)}{2} \text{ for } 1 \leq j \leq k-1.$$

Furthermore for $0 \leq j \leq k-1$

$$u(2k, 2k+2j) = 3k^2 + 2k - \frac{j(5j+3)}{2}$$

and

$$u(2k,2k+2j+1)=3k^2+4k-\frac{j(5j+7)}{2}.$$

The generating function is given by

$$\sum_{k\geq 0} v_3(x,k+2)z^k = \frac{1+x^2z-\left(5+x+2x^2\right)x^3z^2+\left(5x^2+x+2\right)x^4z^3-x^7z^4-x^9z^5}{(1-z)(1-xz)^3(1-x^2z)^2(1-x^3z)} \tag{2.55}$$

For $x=1$ this reduces to $\frac{1+4z+z^2}{(1-z)^4}=\sum_{k\geq 0}(k+1)^3z^k$

and for $x=-1$ we get $\frac{1+6z^2+z^4}{(1-z^2)^3}=\sum_{n\geq 0}(2n+1)^2z^{2n}$.

In the general case we were led to

**Conjecture 2.2**

*For each* $j\geq 1$ *there exists a polynomial* $p_j(z)=p_j(x,z)$ *with degree*

$\deg_z p_j(z)=\frac{(j-1)(j+2)}{2}$ *and* $\deg_x p_j(x,z)=\binom{j+2}{3}-1$ *such that*

$$\sum_{k\geq 0} v_j(x,k+2)z^k = \frac{p_j(x,z)}{(1-z)\prod_{\ell=1}^{j}(1-x^\ell z)^{j+1-\ell}}. \tag{2.56}$$

Moreover it seems that the polynomials have a certain symmetry:

$$(-1)^{\binom{j}{2}}x^{\binom{j+2}{3}-1}z^{\frac{(j+2)(j-1)}{2}}p_j\left(\frac{1}{x},\frac{1}{z}\right)=p_j(x,z).$$

More generally by choosing $x=1$ we get

$$\sum_{k\geq 0} v_j(1,k+2)z^k=\sum_{k\geq 0}(k+1)^j z^k=\frac{p_j(1,z)}{(1-z)^{1+\binom{j+1}{2}}}. \tag{2.57}$$

Recall (cf. [21]) that the Eulerian numbers $\left\langle {n \atop k} \right\rangle$ defined by the recurrence

$$\left\langle {n \atop k} \right\rangle = (k+1)\left\langle {n-1 \atop k} \right\rangle + (n-k)\left\langle {n-1 \atop k-1} \right\rangle$$

with initial values $\left\langle {0 \atop k} \right\rangle = \left[k=0\right]$ and boundary values $\left\langle {n \atop -1} \right\rangle = 0$ satisfy

$$\sum_{k\geq 0}(k+1)^j z^k = \frac{\sum_{\ell=0}^{j}\left\langle {j \atop \ell} \right\rangle z^\ell}{(1-z)^{j+1}}. \tag{2.58}$$

Comparing (2.57) with (2.58) we see that

$$p_j(1,z) = (1-z)^{\binom{j}{2}}\sum_{\ell=0}^{j}\left\langle {j \atop \ell} \right\rangle z^\ell. \tag{2.59}$$

**Conjecture 2.3**

*Let* $r_j(x)$ *be the polynomial defined in* (2.25). *Then*

$$p_j(x,1) = \left(1-x\right)^{j-1}\prod_{\ell=3}^{j}\left(1-x^\ell\right)^{\max(j+1-\ell,0)} r_{j-1}(x). \tag{2.60}$$

For example

$p_2(x,1) = 1 - x^3 = (1-x)r_1(x) = (1-x)(1+x+x^2)$

and

$p_3(x,1) = 1 + x^2 - 5x^3 + x^4 - x^5 + 5x^6 - x^7 - x^9 = (1-x)^2(1-x^3)(1+2x+4x^2+2x^3+x^4).$

## 3. Final remarks

In this paper we studied the sequences

$$a(n,k,t)=\sum_{j\in\mathbb{Z}}(-1)^j\sum_{\ell\geq|j|}\binom{\left\lfloor\frac{n+(k-2)j}{2}\right\rfloor}{\ell-j}\binom{\left\lfloor\frac{n+1-(k-2)j}{2}\right\rfloor}{\ell+j}t^{\ell}.$$

For $t=1$ the numbers $a(n,k,1)$ count the set $A_{n,k}$ of all lattice paths of length $n$ which start at $(0,0)$, stop on heights $0$ or $-1$ and are contained in the strip $-\left\lfloor\frac{k+1}{2}\right\rfloor\leq y\leq\left\lfloor\frac{k}{2}\right\rfloor$.

Other combinatorial interpretations are given by Proposition 1.1 and Proposition 1.2. It also turned out that $a(n,k,1)$ counts all non-negative lattice paths of length $n$ which start at the origin and remain in the strip $0\leq y\leq k$.

For general $t$ there seems to be no connection with path graphs.

But we can find a weight such that the set of all paths from $(0,0)$ to $(n,-n \bmod 2)$ also has weight $c(n,t)=\sum_{j=0}^{n}\binom{\left\lfloor\frac{n}{2}\right\rfloor}{j}\binom{\left\lfloor\frac{n+1}{2}\right\rfloor}{j}t^j$. Let all steps have weight 1 except the down-steps which end on an even height which have weight $t$. Let $c(n,j,t)$ be the weight of all paths from the origin to $(n,j)$. Then we get

$$c(2n,2k,t)=\sum_{j=0}^{n}\binom{n}{j+k}\binom{n}{j}t^j \tag{3.1}$$

and

$$c(2n+1,2k+1,t)=\sum_{j=0}^{n}\binom{n}{j}\binom{n+1}{j+k+1}t^j \tag{3.2}$$

because

$$\sum_{j=0}^{n}\binom{n-1}{j}\binom{n}{j+k}t^j+t\sum_{j=0}^{n}\binom{n-1}{j}\binom{n}{j+k+1}t^j=\sum_{j=0}^{n}\binom{n-1}{j}\binom{n}{j+k}t^j+\sum_{j=0}^{n}\binom{n-1}{j-1}\binom{n}{j+k}t^j$$
$$=\sum_{j=0}^{n}\binom{n}{j+k}\binom{n}{j}t^j$$

and

$$\sum_{j=0}^{n}\binom{n}{j+k}\binom{n}{j}t^j+\sum_{j=0}^{n}\binom{n}{j+k+1}\binom{n}{j}t^j=\sum_{j=0}^{n}\binom{n}{j}\binom{n+1}{j+k+1}t^j.$$

In analogy to the case of numbers there exist polynomials $c(n,j,k,t)$ which satisfy the above relations together with $c(n,j,t,k)=0$ for $j<-\left\lfloor \frac{k+1}{2} \right\rfloor$ and $j>\left\lfloor \frac{k}{2} \right\rfloor$.

Let $A(2n,k,t)=c(2n,0,k,t)$ and $A(2n+1,k,t)=c(2n+1,-1,k,t)$. Then $A(n,k,1)=a(n,k,1)$, but for $k>3$ these polynomials do not coincide with $a(n,k,t)$.

## 3.1. Some Hankel determinants

For fixed $n$ the polynomials $a(n,k,t)$ converge to $a(n,t)$ if $k$ tends to $\infty$. Therefore the Hankel determinants $d(n,k,t)=\det\left(a(i+j,k,t)\right)_{i,j=0}^{n}$ converge to

$$\det\left(a(i+j,t)\right)_{i,j=0}^{n}=t^{\left\lfloor \frac{(n+1)^2}{4} \right\rfloor}.$$

Consider first the special case $t=1$.

Here we have $\det\left(a(i+j,1)\right)_{i,j=0}^{n}=\det\left(\binom{i+j}{\left\lfloor \frac{i+j}{2} \right\rfloor}\right)_{i,j=0}^{n}=1$

and

$$d(n,2k,1)=d(n,2k+1,1)=1 \text{ for } n\le k$$

and

$$d(n,2k,1)=d(n,2k+1,1)=0 \text{ for } n>k.$$

The first identity follows from the fact that for $i,j\le k$

$$\left(a(i+j,k,t)\right)_{i,j=0}^{k}=\left(a(i+j,t)\right)_{i,j=0}^{k} \text{ since } a(n,k,t)=a(n,t) \text{ for } n\le 2k.$$

The second identity is clear because both sequences $a(n,2k,1)$ and $a(n,2k+1,1)$ satisfy a recursion of degree $k+1$.

For general $t$ the situation turns out to be more interesting.

The same argument as above gives

$$d(n,2k,t)=d(n,2k+1,t)=t^{\left\lfloor \frac{(n+1)^2}{2} \right\rfloor}$$

for $n\le k$ and

$$d(n,2k,t)=d(n,2k+1,t)=0 \text{ for } n\ge 2k.$$

But what can be said for the gaps $n\in[k+1,2k-1]$ ?

Computation for small values shows that $d(k+j,2k,t)$ and $d(k+j,2k+1,t)$ are both multiples of $t^{k+j}\left(t-1\right)^{j^2}$.

Let us first consider the $k-1$ polynomials $w(2k)=\left(\dfrac{d(k+j,2k,t)}{t^{k+j}\left(t-1)^{j^2}\right)}\right)_{j=1}^{k-1}$ for $k\geq 2.$

We get

$$w(4)=\left(1\right),\ w(6)=\left(1+t,1\right),\ w(8)=\left(t\left(1+t+t^2\right),\ 1+4t+t^2,\ 1\right),$$
$$w(10)=\left(t^2\left(1+t\right)\left(1+t^2\right),\ t\left(1+4t+10t^2+4t^3+t^4\right),\ (1+t)\left(1+8t+t^2\right),\ 1\right),\cdots.$$

Some of these terms occur in the sequence

$$\left(a(n,t)\right)_{n\geq 0}=\left(1,\ 1+t,\ 1+2t,\ 1+4t+t^2,\ 1+6t+3t^2,\ (1+t)\left(1+8t+t^2\right),\ \cdots\right).$$

Other terms occur in Hankel determinants of $a(n,t).$ This led to

**Conjecture 3.1**

*Let* $a(n,t)=\sum_{j=0}^{n}\binom{\left\lfloor \frac{n}{2}\right\rfloor}{j}\binom{\left\lfloor \frac{n+1}{2}\right\rfloor}{j}t^j$ *and* $D(m,n,t)=\det\left(a(i+j+m,t)\right)_{i,j=0}^{n}$ *for* $n\geq 0$ *and* $D(m,n,t)=1$ *for* $n<0.$

*The Hankel determinants* $d(n,2k,t)=\det\left(a(i+j,2k,t)\right)_{i,j=0}^{n}$ *satisfy*

$d(n,2k,t)=t^{\left\lfloor \frac{(n+1)^2}{2}\right\rfloor}$ *for* $n\leq k$ *and* $d(n,2k,t)=0$ *for* $n\geq 2k.$

*For* $1\leq r\leq k-1$ *they satisfy*

$$d(k+r,2k,t)=t^{k+r}(t-1)^{r^2}D(2r,k-2-r). \qquad (3.3)$$

Let us now consider the $k-1$ polynomials $w(2k+1)=\left(\frac{d(k+j,2k+1,t)}{t^{k+j}\left(t-1\right)^{j^2}}\right)_{j=1}^{k-1}$ for $k\geq 2.$

We get

$$w(5)=(1),\ \ w(7)=(t,1),\ \ w(9)=\left(t^2\left(1+t\right),\ t\left(1+t\right),\ 1\right),$$
$$w(11)=\left(t^4\left(1+t\right),\ t^2\left(1+t\right)\left(1+3t+t^2\right),\ t\left(1+3t+t^2\right),\ 1\right),$$

$$w(13)=\left(t^6\left(1+t+t^2\right),\ t^4\left(1+t+t^2\right)\left(1+3t+t^2\right),\ t^2\left(1+3t+t^2\right)(1+t)\left(1+5t+t^2\right),\ t(1+t)\left(1+5t+t^2\right),\ 1\right),\cdots.$$

Some of these terms are Narayana polynomials as seen from

$$\left(N_n(t)\right)_{n\geq 0}=\left(1,\ t,\ t(1+t),\ t\left(1+3t+t^2\right),\ t(1+t)\left(1+5t+t^2\right),\ t\left(1+10t+20t^2+10t^3+t^4\right),\cdots\right).$$

This led us to

**Conjecture 3.2**

*Let $N_n(t)$ be a Narayana polynomial and let $D_N(m,n,t)=\det\left(N_{i+j+m}(t)\right)_{i,j=0}^{n}$ for $n\geq 0$ and $D_N(m,n,t)=1$ for $n<0.$*

*The Hankel determinants $d(n,2k+1,t)=\det\left(a(i+j,2k+1,t)\right)_{i,j=0}^{n}$ satisfy*

$$d(n,2k+1,t)=t^{\left\lfloor\frac{(n+1)^2}{2}\right\rfloor}\ \textit{for}\ n\leq k\ \textit{and}\ d(n,2k+1,t)=0\ \textit{for}\ n\geq 2k.$$

*For $1\leq r\leq k-1$ they satisfy*

$$d(k+r,2k+1,t)=t^{k+r}(t-1)^{r^2}D_N\left(r,\left\lfloor\frac{k-2-r}{2}\right\rfloor\right)D_N\left(r+1,\left\lfloor\frac{k-3-r}{2}\right\rfloor\right). \tag{3.4}$$